\documentclass[12pt]{article}

\usepackage[cp866]{inputenc}

\usepackage[T2A]{fontenc}

\usepackage{amsfonts,amsthm,amsmath,amssymb,amscd,graphicx}

\usepackage{epsfig}
\usepackage{graphicx}

\begin{document}

\newtheorem{lm}{Lemma}
\newtheorem{theorem}{Theorem}
\newtheorem{prop}{Proposition}
\newtheorem{df}{Definition}
\newtheorem{remark}{Remark}
\newtheorem{corollary}{Corollary}
\newtheorem{ex}{Example}

\begin{center}{\large\bf Wild pseudohyperbolic attractor in a four-dimensional Lorenz system}
\end{center}
\medskip
\begin{center}{\bf Gonchenko S.V.$^1$, Kazakov A.O.$^{2}$, Turaev D.$^{3,2}$ }
\end{center}
\medskip
\begin{center}{$^1$ Lobachevsky State University of Nizhny Novgorod, Russia}
\end{center}
\begin{center}{$^2$ National Research University Higher School of Economics, Nizhny Novgorod, Russia}
\end{center}
\begin{center}{$^3$ Imperial College, London, UK}
\end{center}


\begin{abstract}
We present an example of a new strange attractor which, as we show, belongs to a class of wild pseudohyperbolic spiral attractors. We find this attractor in a four-dimensional system of differential equations which can be represented as an extension of the Lorenz system.
\end{abstract}


\section*{Introduction}

In this paper we build an example of a new strange attractor. We show that it belongs to a class of {\em wild pseudohyperbolic spiral} attractors. A theory of pseudohyperbolic spiral attractors was proposed in \cite{TS98}, however examples of concrete systems of differential equations with such attractors were not known. We perform a series of numerical experiments with the strange attractor which exists in a four-dimensional extension of the classical Lorenz system, and demonstrate that this attractor is indeed pseudohyperbolic, spiral (contains a saddle-focus equilibrium), and wild (contains a hyperbolic set with homoclinic tangencies).
We also discuss the notion of pseudohyperbolicity, as a key property that
ensures the robustness of chaotic dynamics, free from stability windows,
and propose an effective method of numerical verification of the pseudohyperbolicity. The pseudohyperbolicity is a generalization of the hyperbolicity property, which imposes much less restrictions on the system but still guarantees that every orbit in the attractor has maximal positive Lyapunov exponent, both for the system itself and for every close system.

\begin{figure}[tb]
\center{\includegraphics[width=1.0\linewidth]{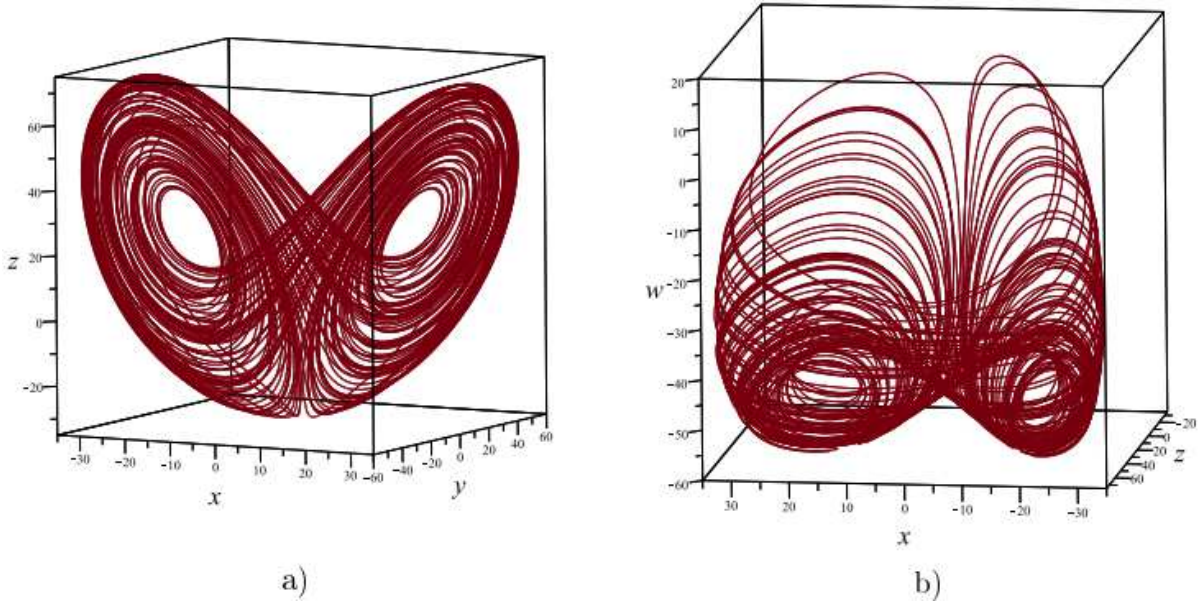} }
\caption{{\footnotesize Projections of the strange attractor existing in system \eqref{eq:LorenzModified} at $\sigma = 10, b = 8/3, r = 25$ and $\mu = 7$ onto: (a) the $(x,y,z)$-plane and (b) the $(x,z,w)$-plane.
}}
\label{fig:WildSpiral}
\end{figure}

We consider the following system of differential equations
\begin{equation}
\left\{
\begin{array}{l}
\dot x = \sigma (y - x), \\
\dot y = x (r-z) - y, \\
\dot z = xy - bz + \mu w, \\
\dot w = -b w - \mu z,
\end{array}
\right.
\label{eq:LorenzModified}
\end{equation}
where $\sigma,r,b$ and $\mu$ are parameters. This system can be viewed as a four-dimensional extension of the classical Lorenz model: when $\mu = 0$ the hyperplane $w=0$ is invariant and, in restriction onto this hyperplane, the system is exactly the Lorenz model. Model \eqref{eq:LorenzModified} was proposed in \cite{book2} (see Part 2, Appendix C, problem C.7.No.86) as a possible candidate for a system with a wild spiral attractor. We perform a series of numerical experiments with the strange attractor which exists in the system at $\mu = 7, \sigma = 10, b = 8/3, r = 25$, see Fig.~\ref{fig:WildSpiral}, and demonstrate that this attractor is indeed pseudohyperbolic and wild.

The pseudohyperbolicity is a key word here. It means that certain conditions hold (see Definition \ref{df:phattr}) which guarantee that every orbit in the attractor is unstable (i.e. it has a positive maximal Lyapunov exponent). Moreover, this instability property persists for all small perturbations of the system.

We recall that one of the main problems of the theory of dynamical systems is that most of the strange attractors discovered in various applications may, in fact, contain stable periodic orbits. These periodic orbits may have quite narrow attraction domains, so we do not see them in numerical experiments, however their existence (either in the system itself or after an arbitrarily small variation of parameters) can be inferred from the existence of homoclinic tangencies \cite{GaS73, N74, G83, GST96, GST08}. In this case one can observe a chaotic behavior (with positive Lyapunov exponent) but can never be sure that increasing the accuracy or the computation time would not make the maximal Lyapunov exponent vanish. Such strange attractors, ``pregnant'' by stable periodic orbits, were called {\it quasiattractors} by Afraimovich and Shilnikov \cite{AfrSh83a}, see also \cite{GST97q}. The corresponding dynamics may appear chaotic for all practical purposes. However, from the purely mathematical point of view, it is a complicated and, quite probably, unsolvable \cite{GST91, GST93a} question whether the dynamics in a given system with a quasiattractor are truly chaotic or become periodic after a long transient process.

Examples of the Afraimovich-Shilnikov quasiattractors are ubiquitous. They include ``torus-chaos'' attractors arising after the breakdown of two-dimensional tori \cite{AfrSh83b} and after a period-doubling cascade, the H\'enon attractor  \cite{Henon76,BenCar}, attractors in periodically perturbed two-dimensional systems \cite{GSV13}, attractors in the Lorenz model beyond the boundary of the region of Lorenz attractor existence \cite{Sh80,BASh92}, spiral attractors in three-dimensional systems with a Shilnikov loop \cite{ACT81,ACT82,ACT85,Chua1986,Gaspard84}, etc. In all these cases, we observe chaotic dynamics but there is no proper mathematical theory which would describe the main properties of such dynamics independently of small perturbations of the system.

However, there exist certain classes of {\it genuinely chaotic attractors} which are not destroyed by small perturbations. These are uniformly hyperbolic attractors, see e.g. book \cite{Kuz_book} and references therein, and Lorenz-like attractors \cite{Lor63,ABS77,Sh81,ABS82,Tucker99, OT17}. Both uniformly hyperbolic and Lorenz attractors are partial cases of {\it pseudohyperbolic} attractors, as proposed in \cite{TS98}.

The following definition generalizes the corresponding definition from \cite{TS98}.
\begin{df} Let a compact set ${\cal A}$ be forward invariant with respect to an $n$-dimensional $C^r$-flow $F$ (i.e.,
$F_t({\cal A}) \subset {{\cal A}}$ for $t > 0$). The set ${\cal A}$ is called \textsf{pseudohyperbolic} if it possesses the following properties.
\begin{itemize}
\item[{\rm 1)}]
For each point $x$ of ${\cal A}$ there exist two continuously dependent on $x$ linear subspaces, $E_1(x)$ with $\dim E_1 = k$ and $E_2(x)$ with $\dim E_2 = n-k$, which are invariant with respect to the differential $DF$ of the flow:
$$
DF_t E_1(x) = E_1(F_t(x)), \qquad DF_t E_2(x) = E_2(F_t(x)),
$$
for all $t \geq 0$ and all $x \in {\cal A}$.
\item[{\rm 2)}] The splitting to $E_1$ and $E_2$ is \textsf{dominated}, i.e., there exist constants $C > 0$ and $\beta > 0$ such that
$$
\|DF_t(x)|_{E_2}\| \cdot \|(DF_t(x)|_{E_1})^{-1}\| \leq C e^{-\beta t}
$$
for all $t \geq 0$ and all $x \in {\cal A}$. (This means that if we have a contraction in $E_2(x)$, then any possible contraction in $E_1(x)$ is uniformly weaker than the contraction in $E_2(x)$, and if we have an expansion in $E_1(x)$, then it is uniformly stronger than any possible expansion in $E_2(x)$).
\item[{\rm 3)}] The linearized flow $DF$ restricted to $E_1$ stretches $k$-dimensional volumes exponentially, i.e., there exist constants $C > 0$ and $\sigma > 0$ such that
$$
det(DF_t(x)|_{E_1}) \geq C e^{\sigma t}
$$
for all $t \geq 0$ and all $x \in {\cal A}$.
\end{itemize}
\label{df:phattr}
\end{df}

Similar definition can be given for diffeomorphisms. Just let the time variable $t$ take discrete values, i.e., $t \in \mathbb{Z}$, and replace $F_t$ in the above definition by the $t$-th iteration of a diffeomorphism $f$, i.e., $F_t = f^t$.

In this paper we consider the case where there is a uniform contraction along the subspaces $E_2(x)$. Thus, using the standard notations of the normal hyperbolicity theory, we will call $E_2(x)$ the strong-stable subspaces and denote them $E^{ss}(x)$; the center-unstable subspaces  $E_1(x)$ will be denoted as $E^{cu}(x)$.

If the pseudohyperbolic set ${\cal A}$ is an attractor, we call it a {\em pseudohyperbolic attractor}. There can be different definitions of an attractor \cite{GorodetskiIlyashenko} but we expect that in reasonable cases the attractor should have an absorbing domain, i.e., a strictly forward-invariant open region ${\cal D}$ that contains ${\cal A}$. We use Ruelle's notion of an attractor \cite{Ruelle81}. Namely, following \cite{GT17} we define the Conley-Ruelle-Hurley (CRH) attractor as a {\it chain-transitive compact invariant set, stable with respect to permanently acting perturbations}\footnote{Recall that the set is called chain-transitive if for any two points in this set and for any $\varepsilon > 0$ there exists an $\varepsilon$-orbit which connects these points. A set is called stable with respect to permanently acting perturbations if for any $\delta > 0$ there exists $\varepsilon > 0$ such that $\varepsilon$-orbits starting at this set never leave its $\delta$-neighborhood.}. Such attractor is always an intersection of a countable sequence of nested absorbing domains.

If all forward orbits from a bounded absorbing domain ${\cal D}$ enter a sufficiently small neighborhood of a pseudohyperbolic attractor ${\cal A}$, then it can be shown that the closure of ${\cal D}$ is also a pseudohyperbolic set. In this case, Condition 3 in Definition 1 obviously guarantees that for every orbit from ${\cal D}$ the maximal Lyapunov exponent is positive. Importantly this property is preserved after any $C^1$-small perturbation. Indeed, since ${\cal D}$ is strictly forward-invariant, it will remain forward invariant for any perturbed system. The dominated splitting Conditions 1--2 are also known to survive \cite{Mane78, Pujals07, PujalsSambarino2009} and the same is obviously true for the volume-expansion Condition 3. Thus, $cl({\cal D})$ remains a pseudohyperbolic set and, even if the attractor $\cal A$ inside ${\cal D}$ changes drastically, it will anyway remain pseudohyperbolic and every orbit of ${\cal A}$ will have positive maximal Lyapunov exponent. In other words, if an attractor is pseudohyperbolic, then stability windows typical for Afraimovich-Shilnikov quasiattractors cannot arise.

In fact, we believe that in the case of diffeomorphisms the following conjecture is true:\\

{\bf P or Q conjecture.} {\it If an attractor is not pseudohyperbolic, it is a quasiattractor.}\footnote{This formulation is very wide. For example, if the attractor is just a stable periodic orbit, it is, formally, a quasiattractor in the sense of our definition. This conjecture becomes meaningful when we speak about attractors of systems with chaotic behavior of unknown nature (when we have a complete knowledge of the structure of the attractor and of its bifurcations, it really does not matter how we name it).}\\

The rationale behind this conjecture is as follows. If we have a chaotic attractor, then it is natural to expect that the attractor should have saddle periodic orbits inside. If the attractor is not pseudohyperbolic, then it is not hyperbolic by definition. Now, in the absence of uniform hyperbolicity one can expect that nontransverse intersections of stable and unstable manifolds of the saddle periodic orbits can be created by small perturbations of the system. It is natural to assume that if the attractor is not pseudohyperbolic, then at least some of such newly created homoclinic tangencies are not pseudohyperbolic\footnote{The closure of a homoclinic orbit is formed by two orbits, the homoclinic orbit itself and the saddle periodic orbit to which it tends both in forward and backward time. These two orbits form a compact invariant set which can be pseudohyperbolic or not according to Definition \ref{df:phattr}.}. In all known cases bifurcations of non-pseudohyperbolic homoclinic tangencies of a diffeomorphism lead to creation of stable periodic orbits \cite{Tat01,GGTat07,GOT14, Gourmelon14}.

It is absolutely not clear how to transform the above arguments to a mathematical proof. Moreover, for the case of flows, there can be mechanisms of pseudohyperbolicity violation other than homoclinic tangencies (e.g. Bykov cycles \cite{Bykov77, Bykov88, Bykov93}) and formulating  the analogous conjecture for flows requires a certain modification of the notion of pseudohyperbolicity \cite{BonattiLuz17}. In any case, it is plausible that without the pseudohyperbolicity or some its extended version (for flows) any chaotic attractor is an Afraimovich-Shilnikov quasiattractor.

In accordance with this philosophy, in order to reliably establish the robustly chaotic dynamics by numerical experiments with a given system, it is not enough to evaluate Lyapunov exponents -- one also needs to check the pseudohyperbolicity of the numerically observed attractor. In term of numerical simulations, if we take a representative trajectory in the attractor and compute Lyapunov exponents $\Lambda_1 \geq \Lambda_2 \geq \dots \geq \Lambda_n$, then Condition 3 from Definition \ref{df:phattr} transforms into
\begin{equation}\label{c1}
\Lambda_1 + \dots + \Lambda_k > 0,
\end{equation}
and Condition 2 becomes
\begin{equation}\label{c2}
\Lambda_k > \Lambda_{k+1}.
\end{equation}
To satisfy the remaining Condition 1 one needs to check that the splitting into a pair of invariant subspaces depends {\em continuously} on the point in the attractor. This requires the computation and analysis of the invariant subspaces $E_1$ and $E_2$ corresponding to the Lyapunov exponents $\Lambda_1, \dots, \Lambda_k$ and $\Lambda_{k+1}, \dots, \Lambda_n$, respectively.

In this paper we propose an effective method of verifying Condition 1, see the description of the method in Sec. \ref{sec:lorph} and test examples in Sec.~\ref{sec:test_examples}. We apply this methodology to system \eqref{eq:LorenzModified}. We show numerically that at $\mu = 7, \sigma = 10, b = 8/3, r = 25$, the system has an absorbing domain with a pseudohyperbolic attractor (with $dim(E^{ss}) = 1$ and $dim(E^{cu}) = 3$).  We also check (see Sec.~\ref{sec:geom}) that the system has a 3-dimensional cross-section in the absorbing domain and the structure of the Poincar\'e map is in agreement with the geometrical model described in \cite{TS98}. Moreover, we verify that the attractor  contains the equilibrium state at zero.

This equilibrium state is a saddle-focus with 1-dimensional unstable manifold and 3-dimensional stable manifold. The fact that this equilibrium is a saddle-focus means that the eigenvalues nearest to the imaginary axis are complex. This implies that the trajectories in the attractor that pass near the saddle-focus have a characteristic spiral shape.

Many examples of strange  attractors where trajectories spiral around a saddle-focus equilibrium have been observed in models of different nature, e.g. R\"ossler system \cite{Rossler}, Arneodo-Coullet-Spiegel-Tresser systems \cite{ACT81,ACT82,ACT85}, Rosenzweig-MacArthur system \cite{KuzFeoRin2001, BakKazKorLevOs2018}, chemical oscillator systems \cite{Gaspard83}, Chua circuit \cite{Chua1986}, etc. The chaoticity of such attractors is explained by the classical Shilnikov theorem \cite{Sh65,Sh70}: if a system has a homoclinic loop to a hyperbolic equilibrium state for which the two nearest to the imaginary axis eigenvalues are complex, then there exists a hyperbolic set in any neighborhood of the homoclinic loop\footnote{This formulation is correct for three-dimensional systems; in higher dimensions one needs additional conditions of general position, see e.g.~\cite{book2}.}. Thus, if we observe a ``spiral attractor'', then we can expect the existence of Shilnikov loop for nearby values of parameters, and the hyperbolic set predicted by Shilnikov theorem can be a part of the attractor.

However, in many cases the spiral attractor is a quasiattractor. For example, in three-dimensional systems of differential equations for which the divergence of vector field is negative (in particular, for all systems mentioned above) every numerically observed spiral attractor must be a quasiattractor. This just follows from the results of \cite{OvsSh87, OvsSh91} that arbitrarily small perturbations of a three-dimensional system with a homoclinic loop to a saddle-focus with negative divergence give rise to stable periodic orbits which coexist with the Shilnikov hyperbolic set.

As our example of system \eqref{eq:LorenzModified} shows, {\em spiral attractors in dimension 4 and higher can carry a pseudohyperbolic structure} and, therefore, be not quasiattractors. Homoclinic loops (and Shilnikov sets) can still be a part of the attractor but the pseudohyperbolicity prevents the birth of stable periodic orbits from such loops.

We believe that in system \eqref{eq:LorenzModified} {\em parameter values corresponding to homoclinic loops to a saddle-focus (Shilnikov loops) are dense in the region of existence of the pseudohyperbolic spiral attractor}, see more discussion on such conjecture in \cite{TS98, TS08}\footnote{In $C^1$-topology this result would follow from Hayashi connecting lemma \cite{Hayashi}; a result from \cite{TS08} provides a $C^{1+\varepsilon}$ version.}. We provide a numerical evidence for this in Sec. \ref{sec:kneading}, see the so-called ``kneading diagrams'' in Fig. \ref{fig:Lorenz4DKneadings}. By \cite{OvsSh87, OvsSh91} bifurcations of such homoclinic loops lead to emergence of homoclinic tangencies. In turn, bifurcations of homoclinic tangencies create the so-called wild hyperbolic sets \cite{N79, GST93b, PV94}.\footnote{The notion of a ``wild hyperbolic set'' was introduced by Newhouse \cite{N70, N79}; this is a uniformly hyperbolic invariant set which has a pair of orbits such that the unstable manifold of one orbit has a nontransversal intersection with the stable manifold of the other orbit in the pair and this property is preserved for all $C^2$-small perturbations -- when we perturb the system, the tangency for a given pair of orbits may disappear, but a tangency between the invariant manifold for another pair of orbits inside the wild hyperbolic set appears inevitably.} Moreover, these wild sets may accumulate to the Shilnikov loops. Since our attractor is the set of all points which are attainable from the saddle-focus equilibrium by $\varepsilon$-orbits for all arbitrarily small $\varepsilon > 0$, it follows that a wild hyperbolic set belongs to the attractor in this case.
Then, the entire unstable manifold of the wild hyperbolic set is also attainable from the saddle-focus and, hence, belongs to the attractor. In particular, the orbits of tangency between the unstable and stable manifolds of the wild hyperbolic set also belong to the attractor\footnote{Under additional assumptions, one can also show that the attractor contains heterodimensional cycles involving saddle periodic orbits with different dimensions of the unstable manifold \cite{Li16,LiT17}. This is a hallmark of the so-called hyperchaos, see e.g. \cite{Rossler79, BK90, Kap95, Stan18, Stan19, GSKK19}.}. It is important because bifurcations of any homoclinic tangency create homoclinic tangencies of arbitrarily high orders, i.e., they cannot be completely described within any finite-parameter unfolding \cite{GST91, GST07}.

Thus, {\em bifurcations of the pseudohyperbolic attractor in system \eqref{eq:LorenzModified} cannot admit a finite-parameter description}. In particular, there can be no good two-parameter description. Therefore, a two-parameter bifurcation diagram (the ``kneading diagram'' presented in Fig.~\ref{fig:Lorenz4DKneadings}) has a characteristically irregular structure. We borrowed the idea of constructing the kneading diagram from \cite{AShil2012, AShil2014}. In these papers, kneading diagrams were built for classical 3-dimensional Lorenz and Shimizu-Morioka systems and it was noted that the kneading diagrams in the regions of existence of the Lorenz attractor have a nice foliated structure, while in the parameter regions where the attractor becomes a quasiattractor the kneading diagrams become ``blurred'', thus indicating the emergence of homoclinic tangencies. A similar blurred structure of the kneading diagram obtained for system \eqref{eq:LorenzModified} confirms the wildness of the pseudohyperbolic attractor we have found in this system.

\section{How to verify the pseudohyperbolicity}  \label{sec:lorph}

The property of pseudohyperbolicity can be expressed in the form of explicitly verifiable cone conditions (see e.g. condition (*) in \cite{ABS77} for Lorenz attractors or Lemma 1 in \cite{TS98} and Theorem 5 in \cite{TS08} for a more general case). This, in principle, opens a way for developing interval arithmetics based numerical tools which could be used for a rigorous establishment of the pseudohyperbolicity (hence, robust chaoticity) of some attractors observed in concrete dynamical systems, similarly to Tucker's computer-assisted proof of the chaoticity of the classical Lorenz attractor \cite{Tucker99}. Such computations are bound to be time-consuming, so one also needs easier to implement less rigorous numerical methods for a fast -- and still reliable -- verification of the pseudohyperbolicity.

The approach we have used in recent papers \cite{GGS12,GGKT14,GG16} is based on computing Lyapunov exponents and checking the fulfillment of inequalities \eqref{c1}, \eqref{c2} for open regions of parameter values, by building the so-called modified Lyapunov diagrams. The idea was that the robustness of conditions \eqref{c1}, \eqref{c2} with respect to parameter changes is an indirect indication of pseudohyperbolicity (providing, in fact, Conditions 2 and 3 of Definition \ref{df:phattr}). In this paper we propose a more reliable approach based on a direct verification of Condition 1 of Definition \ref{df:phattr}.

In our computations we take a very long trajectory of a system, remove a sufficiently long initial segment (to get rid of the transient) and presume that the remaining part of the trajectory gives a good approximation of the attractor. Then we compute the Lyapunov exponents for this piece of the trajectory, along with the corresponding covariant Lyapunov vectors, see more about Lyapunov analysis in \cite{Ginelli2007, Wolfe2007, KuptsovParlitz2012,KupKuz2018}. In such approach the existence of the invariant subspaces $E_1(x)$ and $E_2(x)$ is automatic. So, verifying Condition 1 reduces to checking the continuous dependence of $E_1$ and $E_2$ on the point $x$ in the attractor. If $E_1$ and $E_2$ depend continuously on $x$, then the angle between $E_1$ and $E_2$ stays bounded away from zero (by compactness of the attractor). This observation is used in \cite{KupKuz2018} for verifying the pseudohyperbolicity: one concludes pseudohyperbolicity if the angles between $E_1$ and $E_2$ do not get close to zero.

Our method is different. We plot the graph of the distance between $E_2(x)$ and $E_2(y)$ as a function of the distance between $x$ and $y$ for every pair of points in the attractor (i.e., on the piece of the trajectory which we use for the approximation of the attractor). If $dist(E_2(x), E_2(y)) \to 0$  as $dist(x,y) \to 0$, then we conclude that $E_2$ depends on the point continuously. Importantly, we endow the numerically obtained $E_2$ with an orientation, invariant with respect to the linearized flow, so we measure the distance between oriented spaces $E_2(x)$ and $E_2(y)$. Thus, we check more than required by the pseudohyperbolicity condition 1. Namely, we establish the existence and continuity of an \textit{orientable} field of subspaces $E_2(x)$. Such field may not exist for all pseudohyperbolic attractors (for example, for nonorientable Lorenz attractors \cite{ABS82, SST93}). It always exists when the absorbing domain $\cal D$ (to which the pseudohyperbolicity property of the attractor is extended) is simply-connected. But for a general topology of the attractor, the orientation of $E_2$ may switch when continued along a non-retractable loop. This makes our method applicable to a somewhat narrower class of attractors, however it is enough for our purposes, and taking the orientation into account makes the method more sensitive and reliable, as is seen from the examples bellow.

After the continuity of $E_2$  is verified, we also check the continuity of the field of subspaces $E_1$, also endowed with an invariant orientation. If both the fields $E_1(x)$ and $E_2(x)$ are continuous, we conclude the pseudohyperbolicity of the attractor.

In this paper we consider only the cases when the spaces of strong contraction $E_2(x) = E^{ss}(x)$ are one-dimensional and, thus, the subspaces $E_1(x) = E^{cu}(x)$ have codimension 1. Therefore, the continuity of $E^{cu}(x)$ is equivalent to the continuity of the field of normals $N^{cu}(x)$ to the hyperplanes $E^{cu}(x)$. By the definition, $E^{ss}(x)$ and $N^{cu}(x)$ are line fields; introducing an orientation makes them vector fields. We build the vector fields $\vec E^{ss}(x)$ and $\vec N^{cu}(x)$ by the following numerical procedure.

We consider a system of differential equations
\begin{equation}
\dot x = F(x).
\label{eq:vecField}
\end{equation}
Let ${\cal X} = \{x_1,...,x_m\}$ be a numerically obtained sequence of points on a trajectory of this system corresponding to time moments $t_1,...,t_m$. We compute Lyapunov exponents $\Lambda_1, \dots, \Lambda_n$ for this trajectory and check conditions \eqref{c1}, \eqref{c2}, which in our case take the form
\begin{equation}\label{eq:cond}
\Lambda_1 + \dots + \Lambda_{n-1} > 0,
\end{equation}
\begin{equation}\label{eq:cond0}
\Lambda_{n-1} > \Lambda_{n}.
\end{equation}
Next, we take an arbitrary unit vector $u_m$ at the point $x_m$ and define a sequence of unit vectors $u_s$ at the points $x_s$, $s=1,\dots, m$, by the following inductive procedure: if $u_{s}$ is the vector obtained on the $(m-s)$-th iteration, then $u_{s-1}$ is defined as $u_{s-1} = U_{s-1} / \|U_{s-1}\|$, where $U_{s-1}$ is the solution at $t = t_{s-1}$ of the variation equation
\begin{equation}
\dot U = DF(x(t))\; U
\label{eq:varEq}
\end{equation}
with the initial condition $U(t_s) = u_s$; here $DF$ stands for the matrix of derivatives of $F$ and $x(t)$ is the solution of \eqref{eq:vecField} with the initial condition $x(t_s) = x_s$. We emphasize that we solve equations \eqref{eq:vecField}, \eqref{eq:varEq} in backward time (from $t = t_s$ to $t = t_1$). In order to suppress instability in $x$, we use, at every step, the stored value of $x_s$ as the initial condition, precomputed by integration of \eqref{eq:vecField} in forward time. By \eqref{eq:cond0} the sequence of the unit vectors $u_s$ exponentially converges to the covariant Lyapunov vector corresponding to the Lyapunov exponent $\Lambda_n$, for almost every initial conditions $u_m$. Thus, if $m_1$, $m_2$, and $m$ are sufficiently large, then the segment of the orbit $\cal X$ corresponding to $s \in [m_1, m - m_2]$ gives a good approximation to the attractor and the vectors $u_s$ give a good approximation to $\vec E^{ss}(x_s)$.\footnote{Since on the attractor the sum of all Lyapunov exponents cannot be positive, condition \eqref{eq:cond} implies that $\Lambda_n < 0$, hence the corresponding invariant subspace is indeed contracting.}

We use an analogous procedure to construct vectors $\vec N^{cu}(x_s) = w_s$. We start with a unit vector $w_0$ and define, inductively, $w_{s+1} = W_{s+1}/\|W_{s+1}\|$, where $W_{s+1}$ is the solution at $t = t_{s+1}$ of the adjoint  variation equation
\begin{equation}
\dot W = -[DF(x(t))]^{\top}\; W
\label{eq:varEqTr}
\end{equation}
with the initial condition $W(t_s) = w_s$. Obviously, if $u(t)$ is a solution of (\ref{eq:varEq}) and $w(t)$ is a solution of (\ref{eq:varEqTr}), then the inner product $(u(t),w(t))$ stays constant:
$$\frac{d}{dt} (u,w) = (Au,w) - (u,A^{\top}w)=0$$
(where we denote $A(t)=DF(x(t))$).
Therefore, given any codimension-1 subspace orthogonal to $w_0$, the sequence of its iterations by variational equation (\ref{eq:varEq}) will remain to be orthogonal to $W_s$ at $t=t_s$. Since for a typical choice of such subspace its iterations converge exponentially to $E^{cu}$, it follows that $W_s$ gives a good approximation to $\vec N^{cu}(x_s)$ (orthogonal to $E^{cu}$) for all sufficiently large $s$.

The same procedure works for discrete dynamical systems. We consider a diffeomorphism $x \mapsto F(x)$ and take its trajectory $x_1,\dots,x_m$, where $x_{s+1} = F(x_s)$. Then, the vectors $u_s$ and $w_s$ are determined by the rule
$$
u_{s-1} = \frac{DF(x_s)^{-1} u_s}{\|DF(x_s)^{-1} u_s\|}, \qquad w_{s+1} = \frac{[DF(x_s)^{\top}]^{-1} w_s} {\|DF(x_s)^{\top}]^{-1} w_s\|}.
$$
Note, that the attractor of the map $F$ can have orientable fields of subspaces $\vec E^{ss}$ and $\vec N^{cu}$, but the orientation may flip with each iteration of $F$. To avoid problems with that, we can simply remove from the sequence $(x_s, u_s, w_s)$ every second term.

Finally, once the orbit ${x_s}$, $s \in [m_1, m-m_2]$, and the vectors $u_s$ and $w_s$ are computed, we plot the $\vec E^{ss}$- and $\vec N^{cu}$-\textit{continuity diagrams}. These are graphs in the $(\rho, \varphi)$-plane, where for each pair of points $(x_i, x_j)$, $m_1 \leq i < j \leq m-m_2$,\footnote{In the case of discrete dynamical systems (maps) we consider only even indices $i$ and $j$, to avoid possible problems with orientation flipping.} we plot a point whose coordinate $\rho$ equals to the distance between $x_i$ and $x_j$ and the coordinate $\varphi$ equals to the angle between $u_i$ and $u_j$ for the $\vec E^{ss}$-continuity diagram or between $w_i$ and $w_j$ for the $\vec N^{cu}$-continuity diagram.

These diagrams look like clouds of points in the $(\rho, \varphi)$-plane. If both the $\vec E^{ss}$ and $\vec N^{cu}$ clouds touch the axis $\rho = 0$ only at the single point $(\rho, \varphi) = (0,0)$, then we can conclude that vector fields $\vec E^{ss}(x)$ and $\vec N^{cu}(x)$ are continuous and, thus, the attractor is pseudohyperbolic.

On the other hand, if one of the clouds touches the $\varphi$-axis at nonzero $\varphi$ or there is no visible gap between the cloud and the $\varphi$-axis, then, the corresponding field of subspaces is discontinuous (hence the attractor is not pseudohyperbolic) or it is non-orientable. The latter case may happen only when the cloud touches the axis $\varrho = 0$ just at two points $\varphi = 0$ and $\varphi = \pi$; in this case, one needs more analysis in order to decide whether the attractor is pseudohyperbolic or not.

\section{Test examples} \label{sec:test_examples}

Before applying the method to system (\ref{eq:LorenzModified}), we test it on several examples of strange attractors. We try both the well-known classical models (Lorenz system, H\'enon map, Lozi map, Anosov diffeomorphism) and those that entered the nonlinear dynamics relatively recently (three-dimensional H\'enon maps).

\subsection{Two-dimensional maps.}
In the two-dimensional case the pseudohyperbolicity of the attractor is equivalent to uniform hyperbolicity,
so our method should distinguish between the uniformly-hyperbolic attractors and not uniformly-hyperbolic
ones (the latter can include e.g. Benedicks-Carleson non-uniformly hyperbolic attractor \cite{BenCar,MorViana93,WangYoung08}).

First, we consider the {\em two-dimensional H\'enon map}
\begin{equation}
\begin{array}{l}
\bar x = y, \\
\bar y = M - b x - y^2. \\
\end{array}
\label{eq:Henon2D}
\end{equation}

\begin{figure}[tb]
\center{\includegraphics[width=1.0\linewidth]{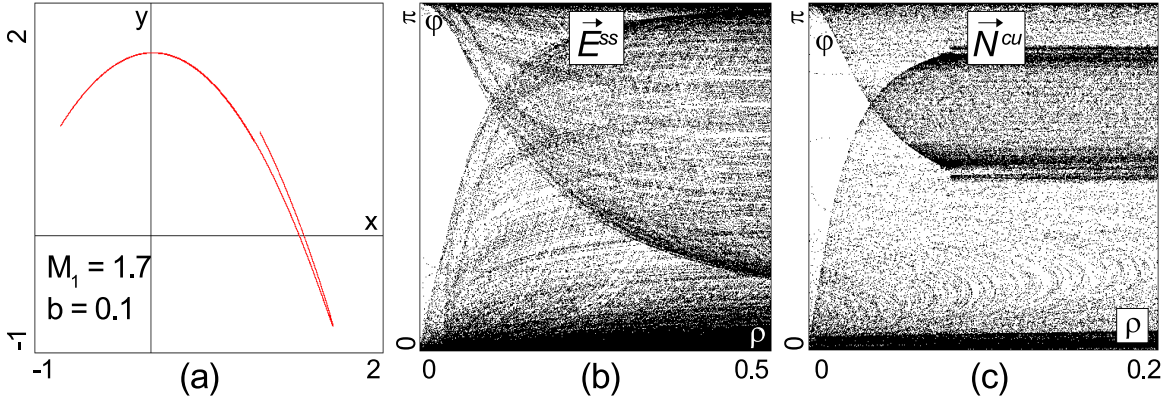} }
\caption{{\footnotesize (a) Attractor of the H\'enon map \eqref{eq:Henon2D} for $b=0.1$, $M = 1.7$. (b) and (c) $\vec E^{ss}$- and $\vec N^{cu}$- continuity diagrams for the attractor.}}
\label{fig:Henon2D}
\end{figure}

In Fig.~\ref{fig:Henon2D} we show numerical results for the H\'enon attractor at
$b=0.1$, $M=1.7$. The attractor is shown in Fig.~\ref{fig:Henon2D}a. The attractor's absorbing domain appears to be simply-connected, so would it be uniformly hyperbolic, the corresponding invariant
line fields $E^{ss}$ and $E^{cu}$ should be orientable. However, since the attractor apparently contains
a saddle fixed point with negative multipliers, each iteration of the map will flip the orientation. Therefore,
in constructing the continuity diagrams we consider only every second iteration of the map.

The continuity diagrams are shown in Figs.~\ref{fig:Henon2D}b and \ref{fig:Henon2D}c. The presence, in these graphs,
of points close to $(0,\varphi)$ with $\varphi$ bounded away from zero indicates the discontinuity of
the vector fields $\vec E^{ss}$ and $\vec N^{cu}$. This confirms the well-known fact that
H\'enon map (as an area-contracting diffeomorphism of a plane) cannot have uniformly-hyperbolic
strange attractors (i.e., any strange attractor in the H\'enon map must be a quasiattractor according to our ``P or Q'' conjecture).

Note that in the $\vec E^{ss}$-continuity diagram the only points close to $\rho=0$
axis are close to $\varphi=0$ or $\varphi=\pi$ -- this means that the line field $E^{ss}$ (i.e., without orientation) would appear continuous here. This demonstrates that taking the orientation of the invariant subspaces into account indeed increases the sensitivity of the method.

\begin{figure}[tb]
\center{\includegraphics[width=1.0\linewidth]{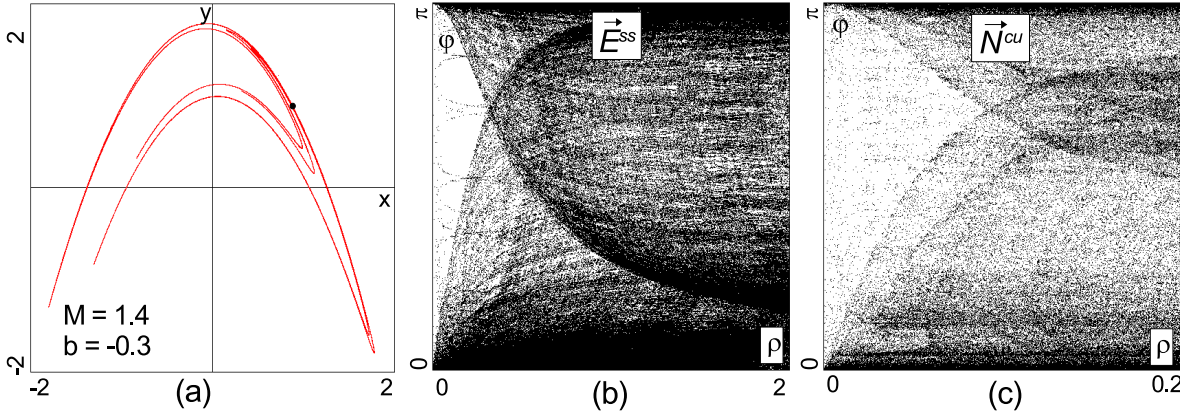} }
\caption{{\footnotesize (a) Attractor of the H\'enon map \eqref{eq:Henon2D} for $b=-0.3$, $M = 1.4$. (b) and (c) $\vec E^{ss}$- and $\vec N^{cu}$- continuity diagrams for this attractor.}}
\label{fig:Henon2D_03}
\end{figure}

In Fig.~\ref{fig:Henon2D_03} analogous results are shown for $b=-0.3$, $M=1.4$.
Both figures~\ref{fig:Henon2D_03}b and \ref{fig:Henon2D_03}c confirm the discontinuity of $E^{ss}$ and $E^{cu}$.

\begin{figure}
\center{\includegraphics[width=1.0\linewidth]{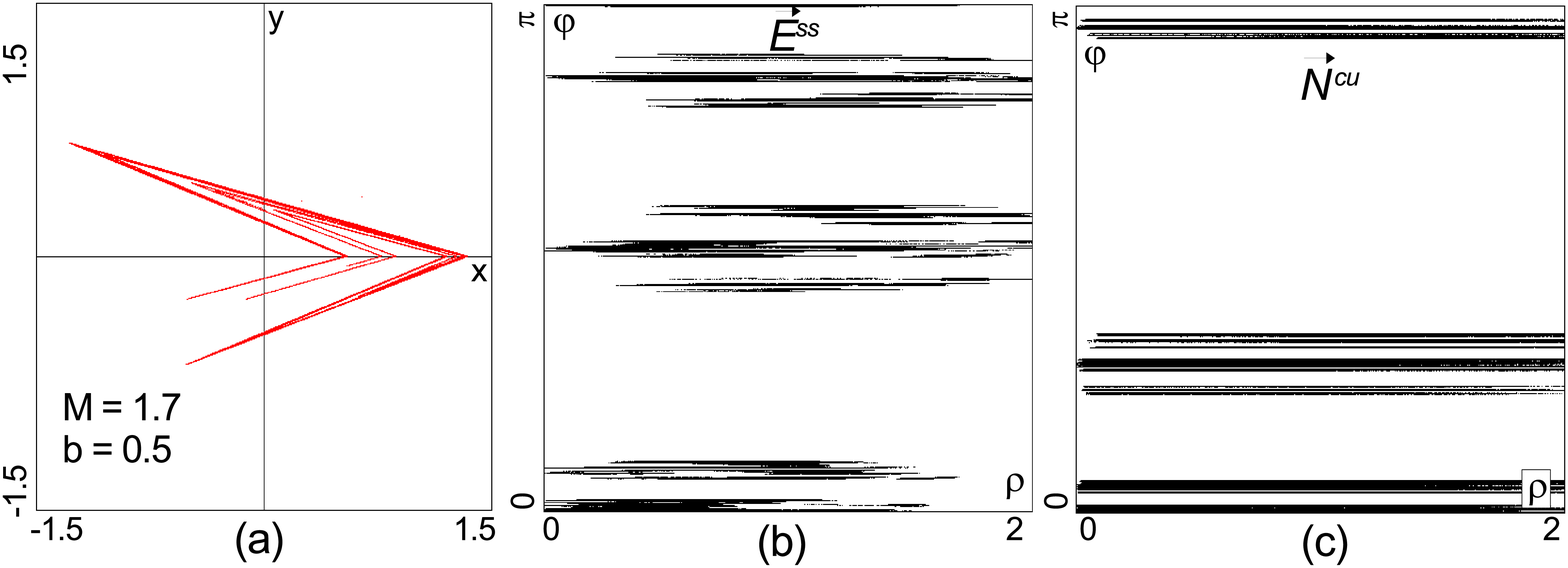} }
\caption{{\footnotesize (a) Lozi attractor of map \eqref{eq:Lozi2D} for $b=0.5$, $M = 1.7$. (b) and (c) $\vec E^{ss}$- and $\vec N^{cu}$- continuity digrams for this attractor.}}
\label{fig:Lozi2D}
\end{figure}

The next example is the {\em Lozi map}
\begin{equation}
\begin{array}{l}
\bar x = 1 +  y - M|x|, \\
\bar y = bx. \\
\end{array}
\label{eq:Lozi2D}
\end{equation}
It is well-known that this map has a singularly-hyperbolic attractor for suitable values of the parameters $M$ and $b$ (e.g. we take $b=0.5$ and $M=1.7$). The singularity appears due to the discontinuity of the derivative at
$x=0$. Thus, we should not expect continuity from the invariant line fields $E^{ss}$ and $E^{cu}$. However, because the map is piecewise affine, the values of the jump in the direction of $\vec E^{ss}$ or $\vec N^{cu}$ at the points of discontinuity must form a certain discrete set.
One can, indeed, clearly see this from Fig.~\ref{fig:Lozi2D} where the $\vec E^{ss}$- and $\vec N^{cu}$- continuity diagrams are formed by horizontal lines that touch the line $\rho=0$ at a certain discrete set of
$\varphi$ values.

Now we consider {\em Anosov diffeomorphisms} of a torus. By definition, these maps are
uniformly hyperbolic. The classical example is
given by the linear map
\begin{equation}
\begin{array}{l}
\bar x = 2 x +  y \;\; (\mbox{mod}\; 1) , \\
\bar y = x + y \;\; (\mbox{mod}\; 1).
\end{array}
\label{eq:An2Tor}
\end{equation}
Both the $\vec E^{ss}$- and $\vec N^{cu}$- continuity diagrams in this case are, quite expectably,
just the lines $\varphi =0$.

\begin{figure}
\center{\includegraphics[width=1.0\linewidth]{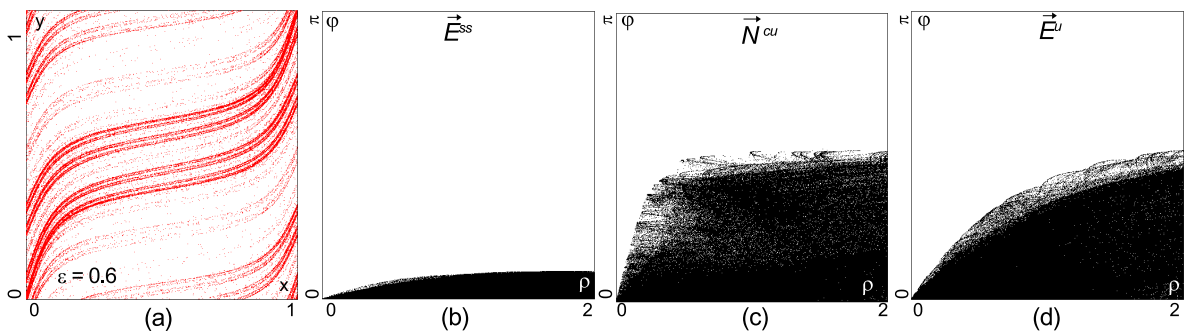} }
\caption{{\footnotesize (a) Attractor of the perturbed Anosov map \eqref{eq:An2Tor_pert} for $\varepsilon=0.6$; (b)--(d) $\vec E^{ss}$- ,  $\vec N^{cu}$- and $\vec E^u$- continuity diagrams. Note, that at $\varepsilon = 0$ all the continuity diagrams flatten just to the straight line $\varphi = 0$.}}
\label{AnosovMobius}
\end{figure}

Small perturbations do not destroy the hyperbolicity of map (\ref{eq:An2Tor}). As an example, we consider
the two-dimensional map from \cite{KaPi19}:
\begin{equation}
\begin{array}{l}
\bar x = 2 \arctan\left(\frac{(1-\varepsilon^2)\sin 2\pi x}{2\varepsilon + (1+\varepsilon^2)\cos 2\pi x} \right) +  y \;\; (\mbox{mod}\; 1) , \\ \\
\bar y = \arctan\left(\frac{(1-\varepsilon^2)\sin 2\pi x}{2\varepsilon + (1+\varepsilon^2)\cos 2\pi x} \right) + y \;\; (\mbox{mod}\; 1).
\end{array}
\label{eq:An2Tor_pert}
\end{equation}
The attractor and the corresponding  $\vec E^{ss}$- and $\vec N^{cu}$- continuity diagrams at $\varepsilon=0.6$ are presented in Fig.~\ref{AnosovMobius}, as well as a similarly constructed continuity diagram for the unstable direction $\vec{E_u}$. The pictures clearly confirm the uniform hyperbolicity of map~\eqref{eq:An2Tor_pert}.

\subsection{Classical Lorenz model.} \label{sec:Lor3D}

An example of a three-dimensional flow which possesses a pseudohyperbolic attractor
for an open set of parameter values is given by
the Lorenz model
\begin{equation}
\left\{
\begin{array}{l}
\dot x = \sigma (y - x), \\
\dot y = x (r-z) - y, \\
\dot z = xy - bz,
\end{array}
\right.
\label{eq:Lorenz}
\end{equation}
where $\sigma$, $r$, and $b$ are parameters. By means of rigorous numerics, it was established by
Tucker \cite{Tucker99} that ``the Lorenz attractor exists'' in this system at
$(\sigma = 10, r = 28, b = 8/3)$. Namely, it follows from the Tucker's result that this system satisfies conditions of the Afraimovich-Bykov-Shilnikov geometrical model \cite{ABS77,ABS82}.

This implies that the attractor of this system at these parameter values is pseudohyperbolic. Thus, there exists a forward invariant absorbing domain ${\cal D}$
within which the Lorenz attractor resides; at each point of $\cal D$ there is a pair of linear spaces $E^{ss}$ and $E^{cu}$ with $\dim E^{ss} = 1$ and $\dim E^{cu} = 2$ such that conditions of Definition~\ref{df:phattr} are satisfied for $E_1 = E^{cu}$ and $E_2 = E^{ss}$.

\begin{figure}[tb]
\center{\includegraphics[width=1.0\linewidth]{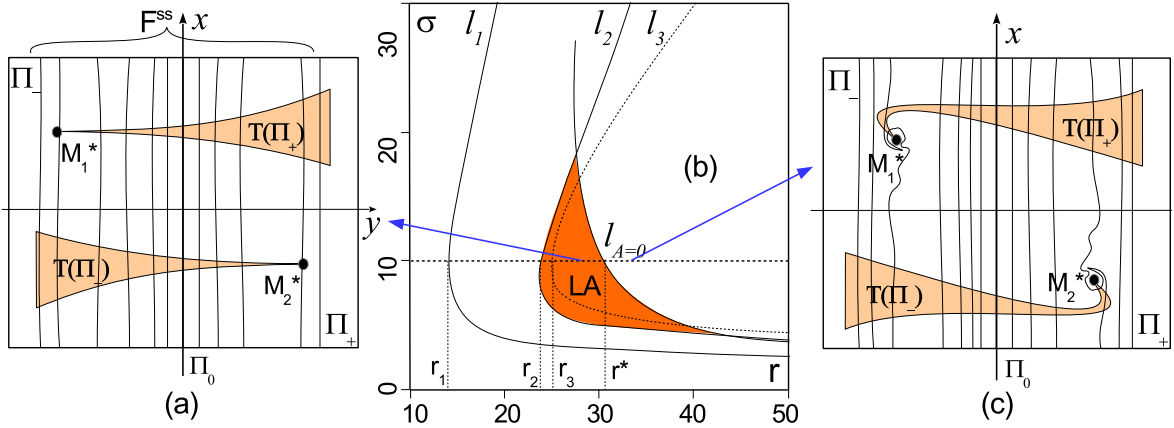} }
\caption{{\footnotesize Poincar\'e map $T$ of the section $\Pi$ (the section $z=r-1$ for the Lorenz model) for values of parameters (a) in $LA$, when the Lorenz attractor exists; (c) to the right of $l_{A=0}$. (b) The domain $LA$ in the $(\sigma,r)$-parameter plane (for $b=8/3$) corresponding to the existence of the pseudohyperbolic Lorenz attractor; the curves $l_1, l_2$ and $l_3$ are described in \cite{Sh80}, the curve $l_{A=0}$ was first computed in \cite{BASh92} and studied in more details in \cite{Creaser17}.}}
\label{fig:foliat}
\end{figure}

By robustness of the pseudohyperbolicity property, the system has the pseudohyperbolic attractor also for some neighborhood of these parameter values. Numerically (non-rigorously) the region $LA$ in the $(\sigma,r)$-parameter plane which corresponds to the existence of the pseudohyperbolic Lorenz attractor for fixed $b=8/3$ was determined in \cite{BASh92, Creaser17}. The left boundary of $LA$ (see Fig.~\ref{fig:foliat}b) is the curve $l_2$ that corresponds to the moment when the unstable separatrices of the saddle equilibrium $O(0,0,0)$ lie on the stable manifolds of certain saddle periodic orbits $L_1$ and $L_2$. These periodic orbits were born from a homoclinic butterfly (to the saddle $O$) which exists when $(\sigma,r)$ belong to the bifurcation curve $l_1$. Along with the orbits $L_{1,2}$ a non-attracting hyperbolic set is born as the homoclinic butterfly splits. This set becomes attracting
(so the Lorenz attractor forms) upon crossing the curve $l_2$ and its attraction basin is bounded by the stable manifolds of $L_1$ and $L_2$.  To the left of $l_2$ the separatrices of $O$ tend to stable equilibrium states $O_1$ and $O_2$, while to the right of $l_2$ they tend to the Lorenz attractor. However, initially, the Lorenz attractor coexists with the stable equilibria $O_1$ and $O_2$; they loose stability on the curve $l_3$
that corresponds to the subcritical Andronov-Hopf bifurcation \cite{Roshchin78}. Here, the periodic orbits $L_1$ and $L_2$ merge with the equilibria $O_1$ and $O_2$ (it happens at $r = r_3\simeq 24.74$ if we fix $b = 8/3$ and $\sigma = 10$). In the region to the right of $l_3$, the equilibria $O_1$ and $O_2$ become saddle-foci with two-dimensional unstable manifolds and the Lorenz attractor becomes the only attractor of the system; see more details in \cite{Sh80, ABS82} and in Chapter 5 of \cite{Shil_book17}.

\begin{figure}
\center{\includegraphics[width=1.0\linewidth]{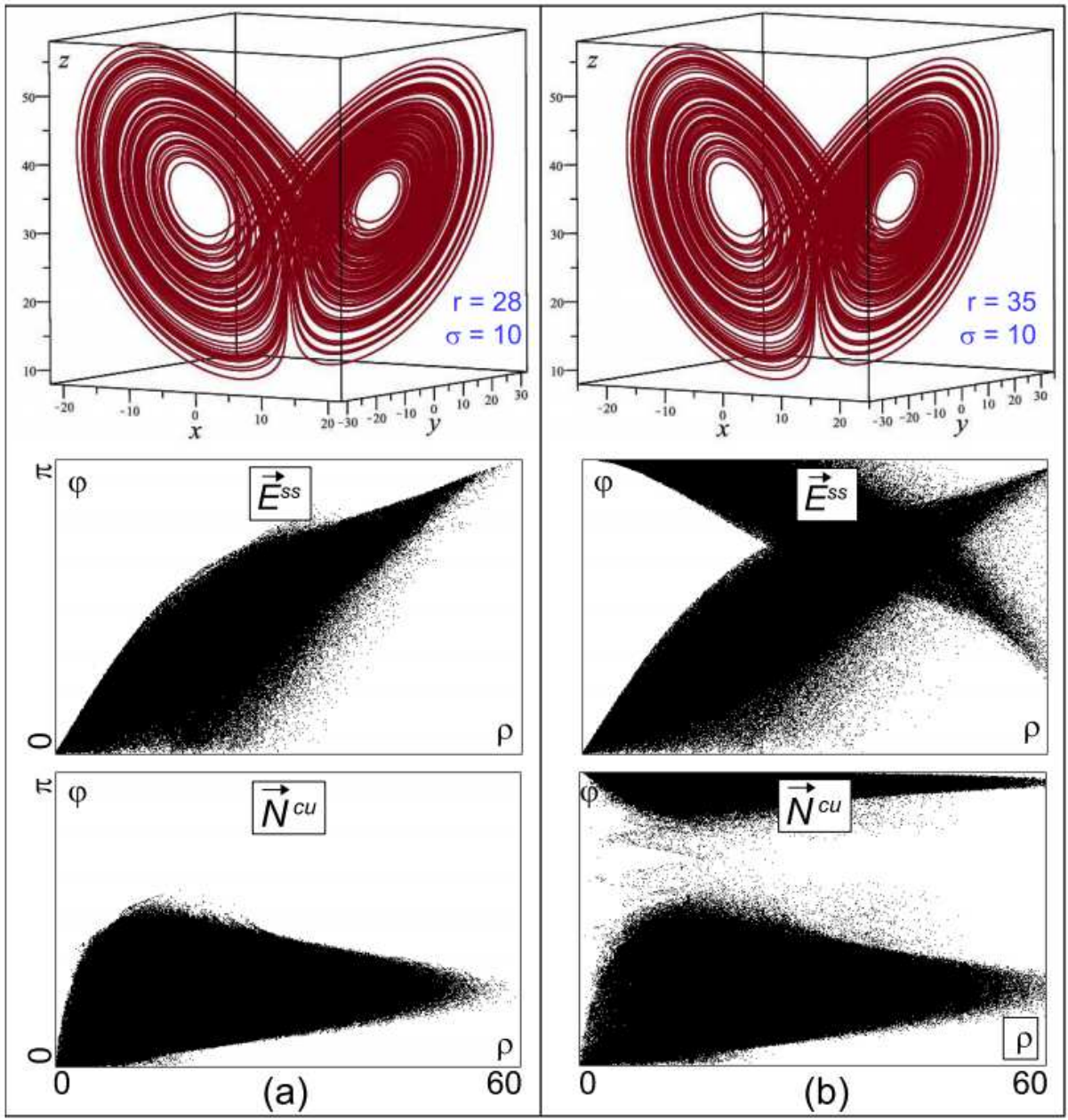} }
\caption{{\footnotesize Different attractors of the Lorenz system (at the top) and the
corresponding continuity diagrams for $\vec E^{ss}(x)$ (in the middle) and for $\vec N^{cu}(x)$ (at the bottom) for the parameter values (a) $r = 28, \sigma  = 10, b = 8/3$, when the attractor is pseudohyperbolic, and (b) $r = 35, \sigma  = 10, b = 8/3$, when Condition 1 of the pseudohyperbolicity (Definition 1) is violated. Note that there is hardly to find any visual difference between the shape of these two attractrors, despite the difference in dynamics.}}
\label{fig:LorenzPseudoTest}
\end{figure}

The right boundary $l_{A=0}$ of the region $LA$ corresponds to the emergence of ``hooks'' in the Poincar\'e map, see Fig.~\ref{fig:foliat}c. System \eqref{eq:Lorenz} has a cross-section, the surface $\Pi:\{z=r-1\}$. The Poincare map $T$ has a discontinuity line $\Pi_0$ corresponding to the intersection of $\Pi$. This line divides the cross-section into two parts, $\Pi_+$ and $\Pi_-$. The images $T(\Pi_+)$ and $T(\Pi_-)$ have a triangular shape, with the vertices at the points $M_-$ and $M_+$ where the unstable separatrices of $O$ intersect $\Pi$ for the first time. Note that the triangles become infinitesimally thin close to the points $M_\pm$. In the region of the existence of the Lorenz attractor, the Poincare map $T$ is (singularly) hyperbolic (the hyperbolicity of the Poincare map is equivalent here to the pseudohyperbolicity of the flow). The hyperbolicity implies the existence of a smooth invariant foliation $F^{ss}$, along which the map $T$ is contracting, see Fig.~\ref{fig:foliat}a. One may conjecture that this foliation still exists at the boundary $l_{A=0}$ and this boundary corresponds to the tangency of the triangles $T(\Pi_\pm)$ at their tip points $M_\pm$ to the foliation. Upon crossing the boundary, the hyperbolicity of the map $T$ gets destroyed. A plausible conjecture is that the curve $l_{A=0}$ is densely filled by points corresponding to the existence of homoclinic loops to $O$ with the so-called separatrix value $A$ equal to zero. Bifurcations of such loops give rise to stable periodic orbits \cite{book2}. Therefore the boundary $l_{A=0}$ separates the region of the pseudohyperbolicity of the Lorenz attractor from the region where it becomes a quasiattractor.

We built $\vec E^{ss}$- and $\vec N^{cu}$- continuity diagrams for the flow of the Lorenz model for parameter values to the left and to the right of the curve $l_{A=0}$, see Fig.~\ref{fig:LorenzPseudoTest}. The diagrams confirm the pseudohyperbolicity of the Lorenz attractor in the region LA and the loss of the pseudohyperbolic structure upon crossing the border
of this region.

Note also that the continuity diagram for $\vec E^{ss}$ touches the line $\rho=0$ only at $\varphi=0$ and $\varphi=\pi$. As we explained, this means either the discontinuity of the line field $E^{ss}$ or its non-orientability. The latter possibility cannot be rejected straight away, as the neighborhood of attractor is not simply-connected (it is a ball around the saddle equilibrium state $O$ and two handles around the two unstable separatrices of $O$). Moreover, bifurcations of a pair of symmetric homoclinic loops with zero separatrix value $A$ can lead to the birth of a non-orientable Lorenz attractor \cite{SST93, Homburg11, BobKazKorSaf19} -- such pairs must exist for parameter values on the boundary curve $L_{A=0}$, so we can predict the existence of ``thin'' non-orientable pseudohyperbolic Lorenz-like attractors  outside the region LA for the parameter values from the so-called ``Shilnikov flames'' \cite{AShil2014}. Nonetheless, the difference between the $\vec E^{ss}$-continuity diagrams in Figs.~\ref{fig:LorenzPseudoTest}a and \ref{fig:LorenzPseudoTest}b clearly indicates that the classical orientable pseudohyperbolic attractor that exists in the region LA is destroyed when the boundary line $l_{A=0}$ is crossed.

\subsection{Lorenz-like attractors in three-dimensional maps.}\label{sec:3DMaps}

\begin{figure}[!th]
\center{\includegraphics[width=1.0\linewidth]{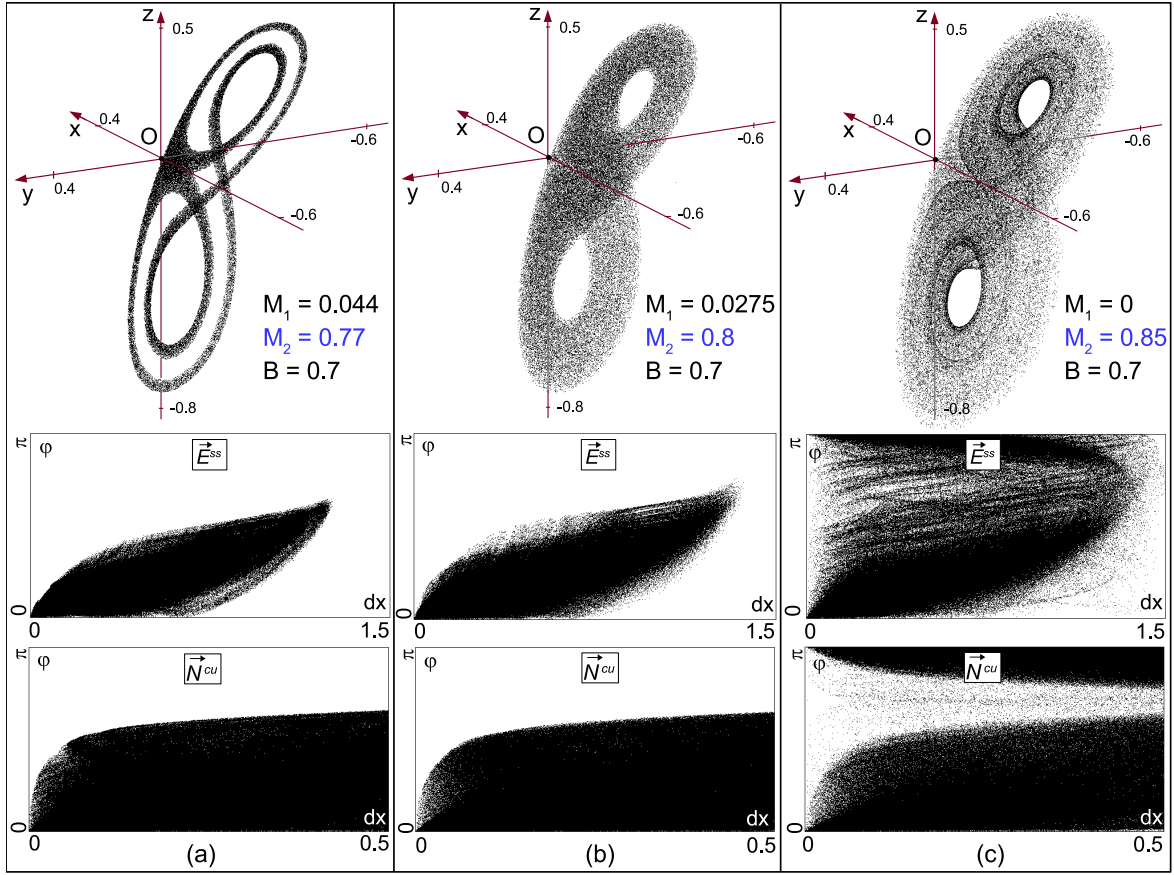} }
\caption{{\footnotesize  Discrete Lorenz-like attractors of map \eqref{3DHM1} (at the top) and their $\vec E^{ss}$-continuity graphs (middle) and
$\vec N^{cu}$-continuity graphs (bottom). Parameter values are (a) $M_1 = 0.044, M_2 = 0.77, B = 0.7$, (b) $M_1 = 0.0275, M_2 = 0.8, B = 0.7$, and (c) $M_1 = 0, M_2 = 0.85, B = 0.7$. Attractor shown in (c) is not pseudohyperbolic.}}
\label{ExDLA1}
\end{figure}

It is shown in \cite{TS08} that adding a small time-periodic perturbation to a system with a pseudohyperbolic attractor does not destroy the pseudohyperbolicity. In particular, the Poincar\'e map (here -- the map over a period of the perturbation) for a small time-periodic perturbation of a system with a Lorenz attractor will
have a {\em discrete Lorenz attractor} -- a pseudohyperbolic attractor similar in shape to the Lorenz attractor of the continuous-time flow \cite{GOST05}. One of the
consequences of this is that discrete Lorenz attractors emerge at local bifurcations of periodic orbits in systems of arbitrary nature. Indeed, it was shown in
\cite{SST93} that the normal form for bifurcations of a periodic orbit with multipliers $(-1,-1,1)$ is a map whose second iteration is the Poincar\'e map
of a small time-periodic perturbation of the Shimizu-Morioka system \cite{SM80}. This system has the (continuous-time) Lorenz attractor for some region of parameter values \cite{ASh86,ASh93}, a rigorous computer assisted proof for this fact was given in \cite{CTZ18}. Therefore, the codimension-3 bifurcation corresponding to a periodic orbit with multipliers $(-1,-1,1)$ can lead, under additional assumptions \cite{SST93,GGOT13}, to the birth of a discrete Lorenz attractor.\footnote{Note that small discrete Lorenz-like attractors can emerge under global bifurcations of multidimensional diffeomorphisms with homoclinic tangencies \cite{GMO06,GOT14} or with nontransversal heteroclinic cycles \cite{GST09,GO13,GO17}. Also universal multi-step bifurcation scenarios leading to discrete Lorenz-like attractors were proposed in the papers \cite{GGS12,GGKT14,GG16}, where their realizations for three-dimensional H\'enon-like maps were also constructed. See also the papers \cite{GGK13,GS19}, where scenarios of the emergence of such attractors were studied in nonholonomic models of Celtic stone.}

An example of such bifurcation was considered in \cite{GOST05} where discrete Lorenz-like attractors were found for the three-dimensional H\'enon map
\begin{equation}
\bar x = y,\;\bar y = z,\; \bar z = M_1 + B x + M_2 y - z^2,
\label{3DHM1}
\end{equation}
in a certain region of the values of parameters $M_1$, $M_2$, and $B$ adjoining to the point $(M_1 = 1/4, M_2 = 1, B = 1)$. This point corresponds to
the existence of a fixed point with multipliers $(-1,-1,1)$, and it was checked in \cite{GOST05} that the normal form for this bifurcation in this map satisfies to the conditions for the birth of the Lorenz attractor. This implies the existence of the pseudohyperbolic attractor for a region of parameter values close enough to this point, see \cite{CTZ18}. However, attractors which look very similar to the Lorenz attractor of the Shimizu-Morioka system were found also
at a sufficient distance from the bifurcation point. For them, the pseudohyperbolicity is not evident and needs to be verified.

In Fig.~\ref{ExDLA1}, examples of discrete Lorenz-like attractors are shown for map \eqref{3DHM1} at $B=0.7$. The continuity diagrams were computed for every second iteration of the map (the map must flips the orientation in $E^{ss}$, as the smallest, i.e., the strongly stable, eigenvalue of the fixed point is negative).
Attractors in Fig.~\ref{ExDLA1}a and ~\ref{ExDLA1}b show the continuity of the field of subspaces $E^{ss}(x)$ and $E^{cu}(x)$, so we can conclude the pseudohyperbolicity, see also \cite{GGKK18} (the necessary conditions $\Lambda_1+\Lambda_2>0$ and $\Lambda_2>\Lambda_3$ were checked in \cite{GOST05}).

In spite of the positivity of the numerically determined in \cite{GOST05} maximal Lyapunov exponent the attractor in Fig.~\ref{ExDLA1}c  is not pseudohyperbolic (the fields of subspaces $E^{ss}(x)$ and $E^{cu}(x)$ are not continuous). In fact, one can show that a stable periodic orbit exists at these parameter value and the ``chaotic attractor'' seen in this Fig.~\ref{ExDLA1}c is an artifact of the (very small) round-off numerical noise \cite{Figueras}.

\section{Pseudohyperbolic spiral attractor} \label{sec:4DLorenz}

The concept of pseudohyperbolic attractors was proposed in \cite{TS98}. In the same paper, a geometric model of the wild spiral attractor for flows in dimension four and higher was constructed.
This geometrical model can be considered as a generalization of the Afraimovich-Bykov-Shilnikov model of the classical Lorenz attractor \cite{ABS77,ABS82}: in the model of \cite{TS98} the saddle equilibrium state of the Lorenz system is replaced by a saddle-focus and the condition of singular hyperbolicity of the Poincar\'e map is replaced by the pseudohyperbolicity.

In \cite{book2}, system \eqref{eq:LorenzModified} was proposed as a possible candidate for a four-dimensional flow which can be described (for some open set of parameter values) by the geometric model from \cite{TS98}. The idea was that at $\mu=0$ system \eqref{eq:LorenzModified} has an invariant three-dimensional hyperplane $w=0$, restricted to which the system is exactly the Lorenz system. So, when we fix the classical Lorenz parameters $r=28, \sigma = 10, b= 8/3$ and take $\mu=0$, system \eqref{eq:LorenzModified} has the Lorenz attractor lying entirely in the hyperplane $w=0$. At small $\mu\neq 0$ the plane $w=0$ is no longer invariant as the saddle equilibrium $O$ at zero becomes a saddle-focus (with a pair of complex conjugate eigenvalues $-b \pm i\mu$), and one numerically observes a strange attractor which includes orbits spiraling around the saddle-focus, see Fig.~\ref{fig:WildSpiral}.

\begin{figure}[tb]
\center{\includegraphics[width=0.7\linewidth]{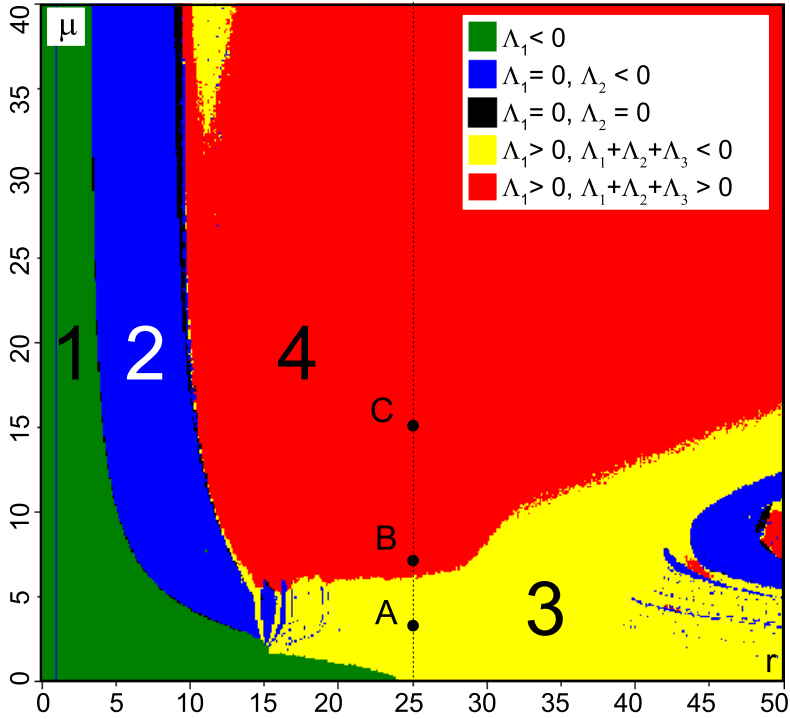} }
\caption{{\footnotesize Diagram of Lyapunov exponents on the $(r,\mu)$-plane for fixed $\sigma = 10, b = 8/3$. Green and blue domains correspond to simple attractors (stable equilibrium and stable limit cycle, respectively). Yellow and red domains correspond to strange attractors with $\Lambda_1+\Lambda_2+\Lambda_3<0$ and  $\Lambda_1+\Lambda_2+\Lambda_3>0$. Note that $\Lambda_1 + \Lambda_2 + \Lambda_3 + \Lambda_4 = -\sigma -2b - 1 < 0$ everywhere in the diagram. }}
\label{Lor4DLyap}
\end{figure}

\begin{figure}[!tb]
\center{\includegraphics[width=0.95\linewidth]{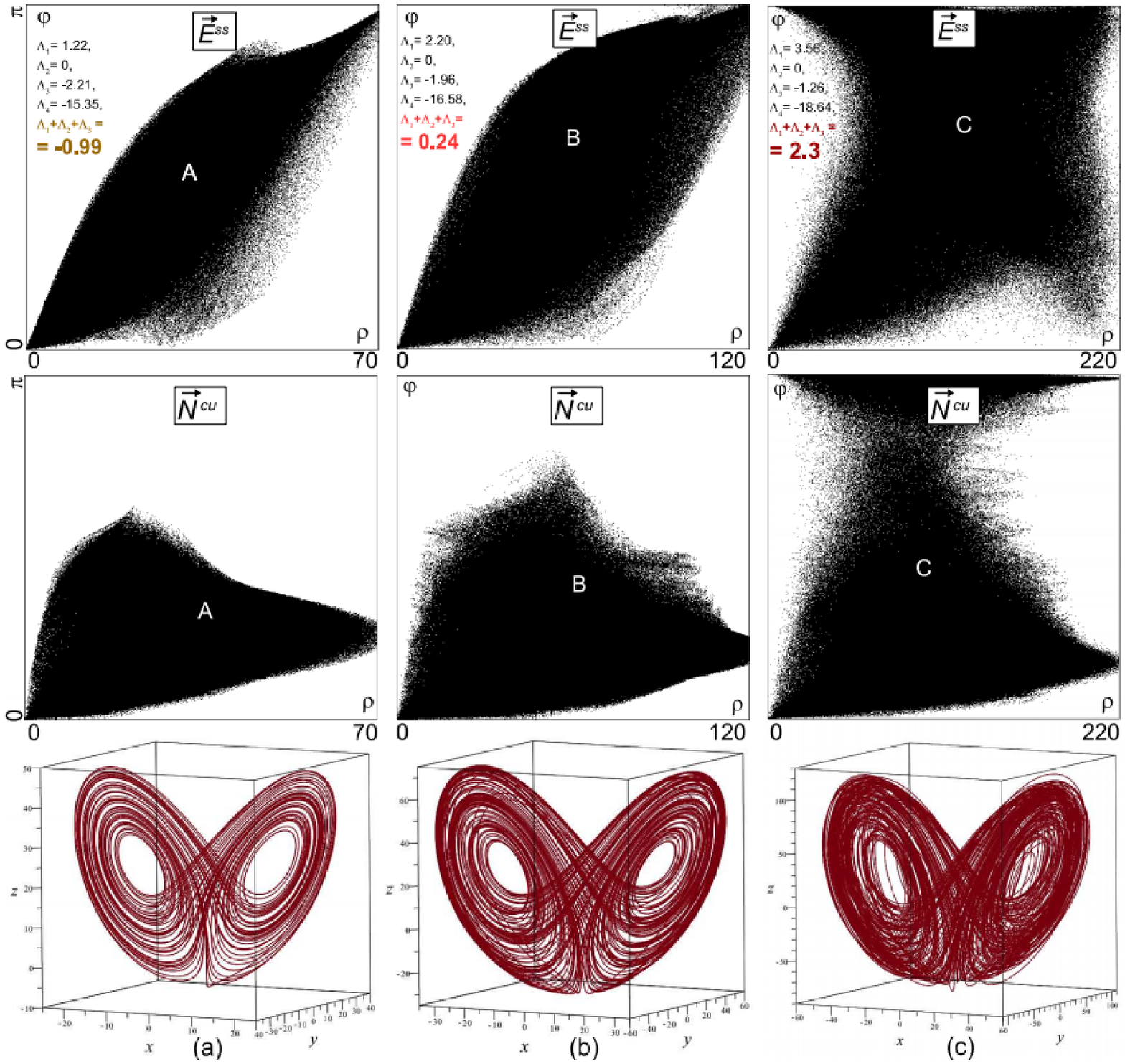} }
\caption{{\footnotesize The $E^{ss}$-continuity diagrams (top), $N^{cu}$-continuity diagrams (middle), and phase portraits of the attractors (bottom) for parameters values corresponding to the points (a) $A\;(r = 25, \mu = 3)$, (b) $B\;(r = 25, \mu =7)$, and (c) $C\;(r = 25, \mu = 15)$ in the $(r,\mu)$-plane at $\sigma = 10$, $b = 8/3$; see the corresponding points in the Lyapunov diagram shown in Fig.~\ref{Lor4DLyap}. The corresponding attractors in bottom (a) and (c) figures, i.e., for parameter points A and C, are quasiattractors (for (a) the necessary condition \eqref{eq:necphd} is not fulfilled, for (c)  $E^{ss}$ and $E^{cu}$ are discontinuous) while the attractor shown in bottom (b) figure is pseudohyprbolic (parameter point B).}}
\label{Lor4DLyap2}
\end{figure}

However, the pseudohyperbolicity conditions are not fulfilled for small $\mu \neq 0$. To see this, note that at $\mu = 0$ the Lorenz attractor for the restriction of the system to the invariant hyperplane $w = 0$ is pseudohyperbolic as guaranteed by the expansion of two-dimensional areas, but at $\mu \neq 0$ we need expansion of three-dimensional volumes. Indeed, the eigenvalues of the linearization at the saddle-focus $O$ are equal to
$$
\begin{array}{l}
\lambda_1 = \frac{1}{2}\left(\sqrt{(\sigma -1)^2 + 4\sigma r} -\sigma -1\right), \\
\lambda_{2,3} = -b \pm i\mu, \\
\lambda_4 = - \frac{1}{2}\left(\sqrt{(\sigma -1)^2 + 4\sigma r} +\sigma +1\right).\end{array}
\label{eq:egval}
$$
At $r=28, \sigma = 10, b= 8/3$ this gives $\lambda_1 \approx 11.83, \lambda_{2,3} = -8/3 \pm i \mu, \lambda_4 \approx -22.83$. Therefore, the space $E^{ss}$ at the point $O$ is one-dimensional (it corresponds to the smallest eigenvalue $\lambda_4$). By continuity of $E^{ss}$, would we have a pseudohyperbolic attractor the space $E^{ss}$ would be one-dimensional at every point of the attractor. Accordingly, the space $E^{cu}$ must be three-dimensional. This condition is not satisfied for small $\mu$ -- the sum of the first three Lyapunov exponents is negative. Indeed, it is well known that the first two Lyapunov exponents for the Lorenz system at the classical parameter values are $\Lambda_1 \approx 0.906$ and $\Lambda_2 = 0$. In system \eqref{eq:LorenzModified} at $\mu = 0$ the Lyapunov exponents $\Lambda_1$ and $\Lambda_2$ remain the same and $\Lambda_3 = -8/3$. This gives $\Lambda_1+ \Lambda_2 + \Lambda_3 \approx -1.761 < 0$ and it cannot become positive for small $\mu$.

In this paper we show that the pseudohyperbolicity is gained for a certain interval of sufficiently large values of $\mu$. We also slightly deviate from the classical value of $r=28$. For example, the pseudohyperbolic attractor is found at $\sigma = 10, b = 8/3, r = 25$,  and $6 < \mu < 12$.

Fig.~\ref{Lor4DLyap} shows a diagram of Lyapunov exponents on the $(r,\mu)$-parameter plane for attractors of system (\ref{eq:LorenzModified}) with fixed $\sigma = 10$ and $b = 8/3$. Different colors correspond to different spectra of the Lyapunov exponents $\Lambda_1> \Lambda_2> \Lambda_3> \Lambda_4$ and, respectively, to different dynamical regimes. Green domain 1 and blue domain 2 correspond to the existence of regular attractors: a stable equilibrium ($\Lambda_1<0$) and a stable limit cycle ($\Lambda_1=0,\Lambda_2<0$), respectively. Yellow domain 3 and red domain 4 correspond to the existence of strange attractors, where $\Lambda_1>0,\Lambda_2=0$ and $\Lambda_1+\Lambda_2+\Lambda_3<0$ in the yellow domain 3 while $\Lambda_1+\Lambda_2+\Lambda_3>0$ for the red domain 4. Note that the numerically observed attractors for the values of parameters from domains 3 and 4 are always spiral attractors that appear to contain the saddle-focus $O$ at zero. The necessary condition for the pseudohyperbolicity of the attractor
\begin{equation}
\Lambda_1>0, \Lambda_2 =0, \Lambda_1+\Lambda_2+\Lambda_3 >0
\label{eq:necphd}
\end{equation}
is fulfilled only in the red domain 4. In particular, the attractor in Fig.~\ref{Lor4DLyap2}a (point A in the diagram of Fig.~\ref{Lor4DLyap}) is not pseudohyperbolic.

The focus of our investigation will be the attractor at
\begin{equation}
r = 25, \sigma = 10, b = 8/3, \mu = 7.
\label{eq:Lor4DParams}
\end{equation}
The corresponding point $(r=25,\mu =7)$ (point B) belongs to domain 4 from Fig.~\ref{Lor4DLyap}. Therefore, the attractor satisfies necessary condition \eqref{eq:necphd} for pseudohyperbolicity (the numerically obtained exponents are $\Lambda_1 \approx 2.19, \Lambda_2 \approx 0,  \Lambda_3 \approx -1.96, \Lambda_4 \approx -16.56$).

The attractor is shown in Fig.~\ref{Lor4DLyap2}b. To establish its pseudohyperbolicity, we need to verify that the subspaces $E^{ss}(x)$ and $E^{cu}(x)$ depend continuously on the point of the attractor. We did it by computing the $\vec E^{ss}$- and $\vec N^{cu}$-continuity diagrams, as discussed in Section \ref{sec:lorph}. The diagrams are shown in Fig.~\ref{Lor4DLyap2}b. They are quite similar to those for the Lorenz attractor (compare Figs.~\ref{Lor4DLyap2}b and \ref{fig:LorenzPseudoTest}a) and clearly show the sought continuity of $E^{ss}$ and $E^{cu}$.

Note that at the further increase of $\mu$  the continuity condition gets broken, i.e., the attractor loses the pseudohyperbolicity. For example, the attractor shown in Fig.~\ref{Lor4DLyap2}c corresponds to point C ($r=25, \mu = 15$) in the diagram from Fig.~\ref{Lor4DLyap}. Here, the $\vec E^{ss}$- and $\vec N^{cu}$- continuity diagrams clearly indicate the lack of continuity.

\subsection{Spiral geometry of the attractor.} \label{sec:geom}

In the rest of the paper we study dynamical properties of the pseudohyperbolic attractor found for the parameter values given by (\ref{eq:Lor4DParams}). First, we establish that the system has an absorbing domain that contains $O$ and has a special structure similar to that described in \cite{TS98}.

\begin{figure}[tb]
\center{\includegraphics[width=1.0\linewidth]{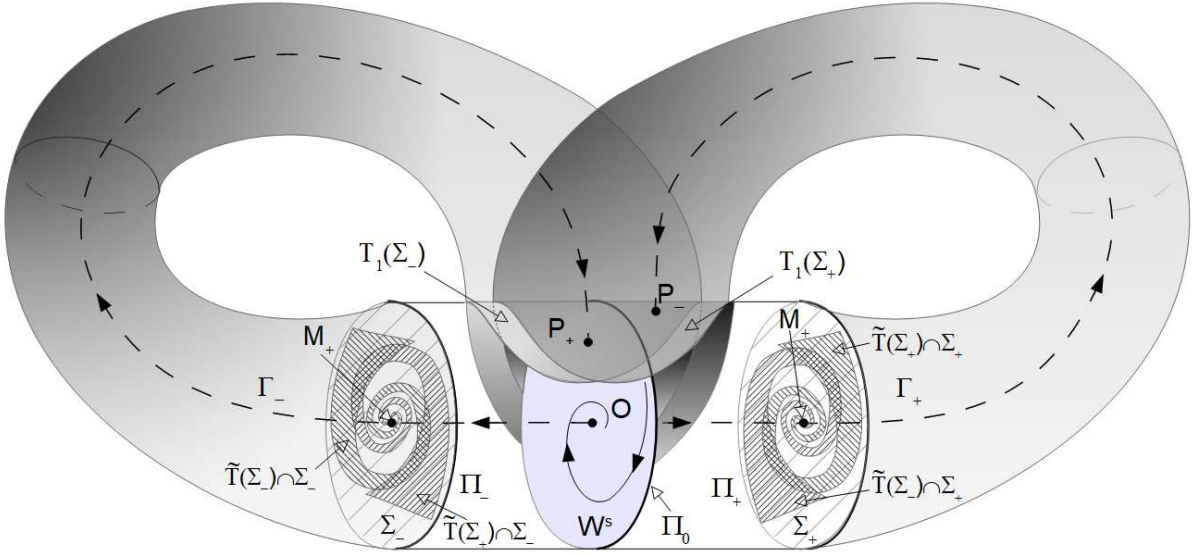}}
\caption{{\footnotesize The scheme of Poincar\'e maps for the wild spiral attractor.}}
\label{fig:PM_cyl}
\end{figure}

In \cite{TS98} the system is assumed to have a cross-section $\Pi$, a three-dimensional cylinder whose intersection with $W^s(O)$ contains a two-dimensional annulus $\Pi_0$ which divides $\Pi$ into two cylinders $\Pi_+$ and $\Pi_-$. Both unstable separatrices $\Gamma_+$ and $\Gamma_-$ of $O$ are assumed to intersect $\Pi$. We denote as $P_+$ and $P_-$ the points of the first intersection of $\Gamma_+$ and, respectively, $\Gamma_-$ with $\Pi$. Thus, the orbits starting in $\Pi$ near
$\Pi_0$ return to $\Pi$ near the points $P_+$ and $P_-$. Moreover, we also assume that all the orbits starting in $\Pi_+ \cup \Pi_-$ return to $\Pi$. Thus, the Poincar\'e map $T: \Pi_{+} \cup \Pi_{-} \rightarrow \Pi$ is defined, see Fig.~\ref{fig:PM_cyl}.
In this construction, if we take the union of all forward orbits starting in $\Pi$ and add to it the two separatrices $\Gamma_+$ and $\Gamma_-$, then we obtain an absorbing domain.

It happened to be difficult to find an explicit expression for the cylindric cross-section $\Pi$ in our system. However, in the above described construction, instead of the cylinder $\Pi$ we may take, as a cross-section, a pair of disjoint balls $\Sigma_+$ and $\Sigma_-$  transverse to $\Gamma_+$ and $\Gamma_-$, respectively, see Fig.~\ref{fig:PM_cyl}. Since every point starting at $\Pi_+$ before returning to $\Pi$ must intersect $\Sigma_+$ and every point starting at $\Pi_-$ before returning to $\Pi$ must intersect $\Sigma_-$, the analysis of the Poincar\'e map on $\Pi$ is equivalent to the analysis of the Poincar\'e map $\tilde T$ on $\Sigma = \Sigma_+ \cup \Sigma_-$.

\begin{figure}[!tb]
\center{\includegraphics[width=1.0\linewidth]{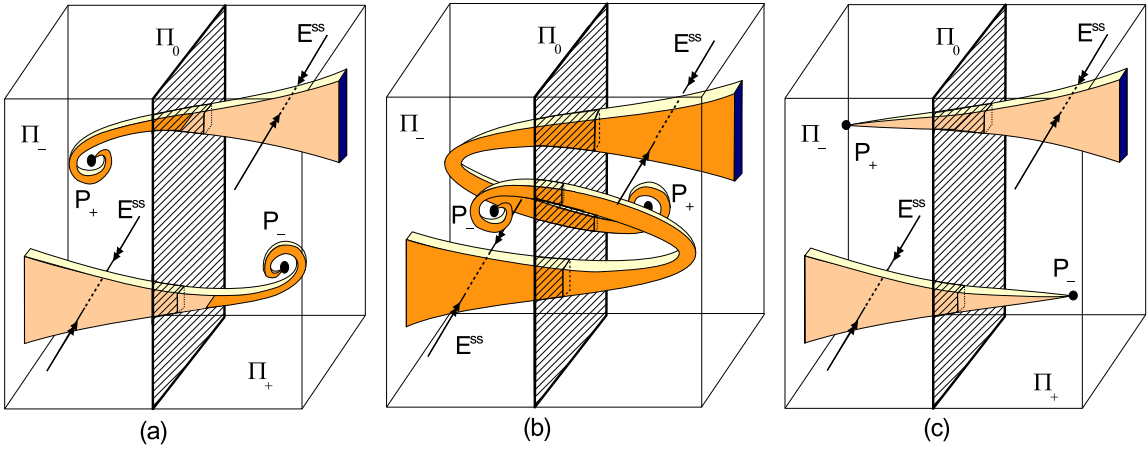} }
\caption{{\footnotesize (a) At small $\mu$ the Poincar\'e map expands two-dimensional areas transverse to the strong stable direction only near $\Pi_0$ (the image of this region is shown by a darker color). The attractor is not pseudohyperbolic (the sum of the three largest Lyapunov exponents is negative). (b) The attractor becomes pseudohyperbolic at larger $\mu$: the Poincar\'e map expands areas transverse to $E^{ss}$ everywhere in a neighborhood of the attractor. (c) Scheme of the Poincar\'e map at $\mu=0$.}}
\label{figfoliat4d}
\end{figure}

We can represent the map $\tilde T$ as $\tilde T= T_0 \circ T_1$, where $T_1$ takes $\Sigma$ into $\Pi$ and $T_0$ takes $\Pi_+$ into $\Sigma_+$ and $\Pi_-$ into $\Sigma_-$. The image $T_1(\Sigma_+)$ in $\Pi$ is divided by $\Pi_0$ into two regions, one further goes to $\Sigma_+$, the other goes to $\Sigma_-$. Since orbits passing near $W^s(O)$ come close to the saddle-focus and, therefore, spiral around the unstable separatricies $\Gamma_+$ and $\Gamma_-$, the image $\tilde T(\Sigma_+) \cap \Sigma_+$ has a form of a wedge spiraling to the point $M_+=T_1^{-1}P_+$ and the image $\tilde T(\Sigma_+) \cap \Sigma_-$ has a form of a wedge spiraling to the point $M_-=T_1^{-1}P_-$. The same is true for the image of $\Sigma_-$. On the cross-section $\Pi$ the images of these spiral wedges by the map $T_1$, i.e., the set $T_1 \circ \tilde T(\Sigma)$, have the form schematically presented in Fig.~\ref{figfoliat4d}a,b. Would the equilibrium $O$ be a saddle instead of the saddle-focus, the picture would be as shown in Fig.~\ref{figfoliat4d}c, i.e., the same as in the Lorenz model with an additional contracting direction (compare with Fig.~\ref{fig:foliat}a).

Thus, Fig.~\ref{figfoliat4d}c depicts the action of the Poincar\'e map for system (\ref{eq:Lor4DParams}) at $\mu=0$, while Fig.~\ref{figfoliat4d}a shows the behavior at small $\mu\neq 0$. As we mentioned, the sum of the three largest Lyapunov exponents at $\mu= 0$ is negative, and it cannot become positive for small $\mu\neq 0$, therefore we do not have pseudohyperbolicity at small $\mu\neq 0$. Namely, the Poincar\'e map does not expand two-dimensional areas transverse to $E^{ss}$ (it can expand the two-dimensional areas only near the surface $\Pi_0$). Our understanding for the onset of pseudohyperbolicity as $\mu$ grows to non-small values is that at such values of $\mu$ the images  $T_1 \left( \tilde T(\Sigma) \cap \Sigma_+\right)$ and $T_1 \left( \tilde T(\Sigma) \cap \Sigma_-\right)$ start spiral around their tips $P_+$ and $P_-$ with higher amplitude, giving enough room for the expansion of areas everywhere on the cross-section, as shown in Fig.~\ref{figfoliat4d}b.

\begin{figure}[tb]
\center{\includegraphics[width=0.9\linewidth]{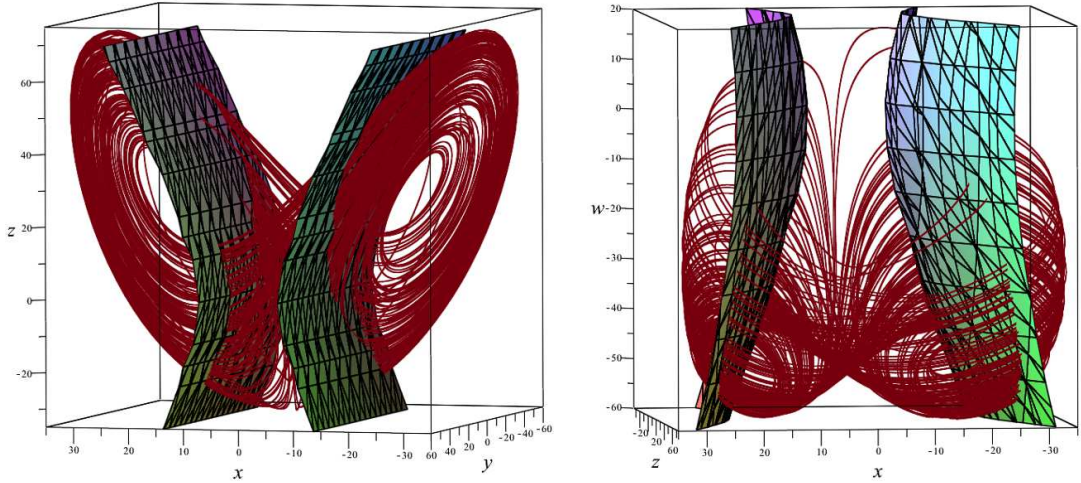} }
\caption{{\footnotesize Projections of the attractor and cross-section $\Sigma$ onto (a) the hyperplane $(x, y, z)$ and (b) the hyperplane $(x, z, w)$.}}
\label{Fig:3DSectionChoice}
\end{figure}

In order to construct the cross-section $\Sigma$ in system \eqref{eq:LorenzModified} we take the three-dimensional hypersurface
\begin{equation}
z = \sqrt{9x^2 - w^2 - 550}.
\label{eq:sec}
\end{equation}
The boxes $\Sigma_+$ and $\Sigma_-$ are the parts of this surface defined by the inequalities
\begin{equation}
x \in [10, 30], y \in [-20, 20], w \in [-60, -10]
\label{eq:cond1}
\end{equation}
and, respectively,
\begin{equation}
x \in [-30, -10], y \in [-20, 20], w \in [-60, -10].
\label{eq:cond2}
\end{equation}

We check that the orbits of the flow intersect such chosen $\Sigma$ transversely, see Fig.~\ref{Fig:3DSectionChoice}. It is worth noting  that we can not use the hyperplane $z = const$ as a cross-section. Such choice, inherited from the Lorenz system, could be natural at small $\mu$, but in the case of non-small $\mu$ we consider here the orbits that wind around the saddle-focus inevitably touch such planes. In general, the problem of choosing a good cross-section in problems of such kind is not trivial.

\begin{figure}[tb]
\center{\includegraphics[width=0.5\linewidth]{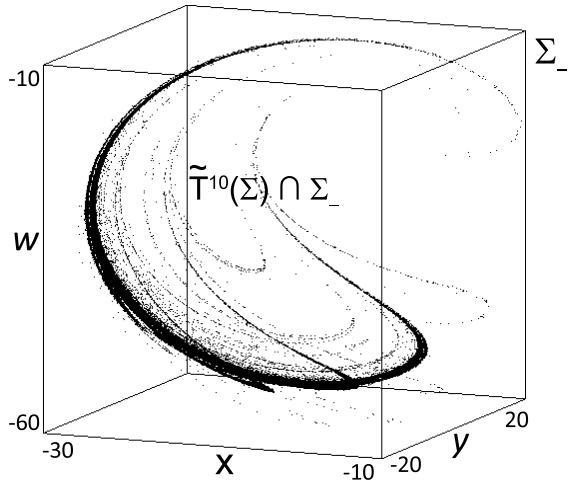} }
\caption{{\footnotesize The image of the cross-section $\Sigma_-$ after $10$ iterations of the return map in the intersection with $x <0$ lies strictly inside $\Sigma_-$.}}
\label{fig:PM_10iter}
\end{figure}

In our case, we encountered a problem that not all the orbits starting in $\Sigma$ return inside it. Namely, all the points starting in $\Sigma$ return to the hypersurface \eqref{eq:sec}, however not all of them satisfy conditions \eqref{eq:cond1} or \eqref{eq:cond2}. We resolve this issue by considering the $10$-th return to the hypersurface \eqref{eq:sec}. For the uniform grid of $200 \times 200 \times 200$ of initial conditions on $\Sigma_+$ and $\Sigma_-$, we checked that the image after the $10$-th return to the hypersurface lies strictly inside $\Sigma$, see Fig.~\ref{fig:PM_10iter}. This confirms the existence of the absorbing domain.

\begin{figure}[th]
\center{\includegraphics[width=1.0\linewidth]{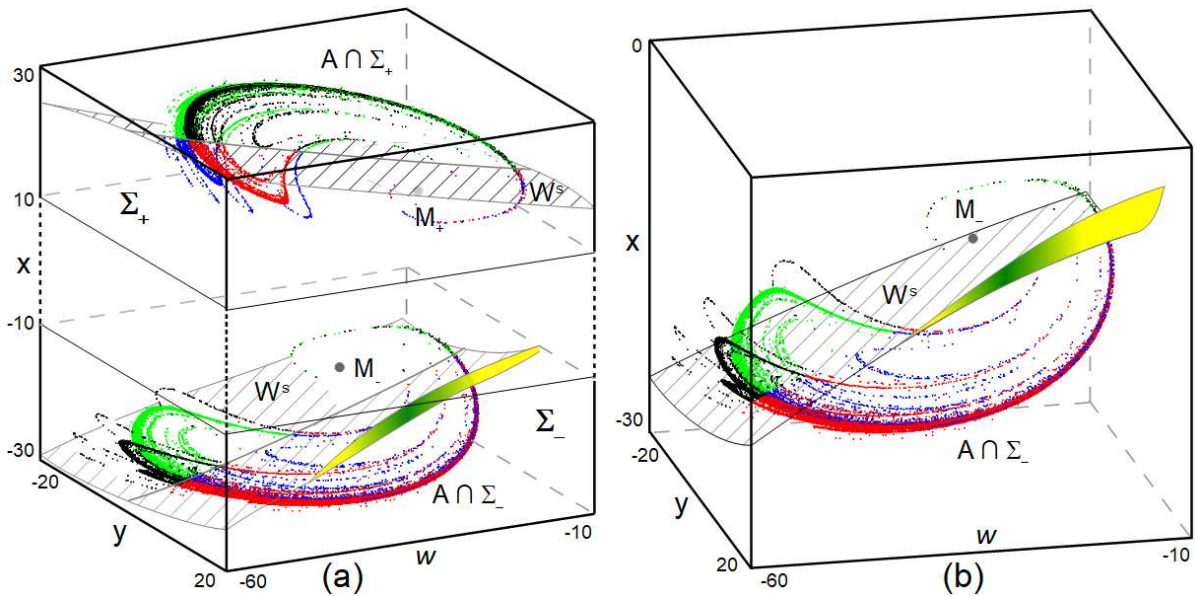} }
\caption{{\footnotesize (a) Attractor of system \eqref{eq:LorenzModified} for parameter values \eqref{eq:Lor4DParams} in the intersection with the cross-section $\Sigma$. Green and black points are those whose iterations by $\tilde T$ belong to $\Sigma_-$; red and blue points are those whose images by $\tilde T$ belong to $\Sigma_+$. The images of green and blue points by $\tilde T^{-1}$ belong to $\Sigma_+$ and the images of red and black points by $\tilde T^{-1}$ belong to $\Sigma_-$. The surface $W^s$ is a piece of the stable manifold of the point $O$ defined as $T_1^{-1}\Pi_0$; it separates green and black poits in the attractor from red and blue ones. The visible presence of the intersection of the attractor with the stable manifold of the saddle-focus $O$ confirms that $O$ belongs to the attractor, i.e., this is a spiral attractor. (b) The part of the attractor that lies in $\Sigma_-$ (an enlarged version of the corresponding fragment of (a)).}}
\label{Fig:3DSection}
\end{figure}

We also checked that the one-dimensional unstable separatrices of $O$ intersect $\Sigma$; the intersection points $M_\pm$ are shown in Fig.~\ref{Fig:3DSection}. In the same figure we show the attractor $A$ of the separatrices, obtained numerically by computing $6\cdot 10^5$ intersections of the separatrices with $\Sigma$ and omitting the first $10^5$ intersection points\footnote{One can think of the numerically obtained attractor as an approximation of the $\omega$-limit set of the separatrices. However, the numerical trajectories are, actually, epsilon-orbits, so it is safer to think of this attractor as a prolongation of the separatrices \cite{GT17}, i.e., the set of points attainable from the saddle-focus by $\varepsilon$-orbits for arbitrarily small $\varepsilon$.}. We use the following color coding: the images of green and black points by $\tilde T$ belong to $\Sigma_-$ and the images of red and blue points by $\tilde T$ belong to $\Sigma_+$ while the images of green and blue points by $\tilde T^{-1}$ belong to $\Sigma_+$ and the images of red and black points by $\tilde T^{-1}$ belong to $\Sigma_-$.

Obviously, the boundary between ``green and black'' and ``red and blue'' points corresponds to the intersection of $\Sigma$ with $T_1^{-1} \Pi_0$, a piece of the stable manifold $W^s(O)$. We computed this surface by a numerical procedure independent of the computation of the attractor. We took the uniform grid of $200\times 200 \times 200$ initial points and interpret as $W^s$ the boundary between the points whose first iteration by the Poincar\'e map $\tilde T$ lies in the region $x<0$ and the points whose first iteration by $\tilde T$ lies in $x>0$.

It is clearly seen in Fig.~\ref{Fig:3DSection} that the attractor $A$ intersects the surface $W^s$. Therefore, we can conclude that the attractor of the separatrices of $O$ for the flow of system \eqref{eq:LorenzModified} for the parameter values given by \eqref{eq:Lor4DParams} intersects the stable manifold of $O$. Hence,
it contains the saddle-focus $O$ itself. This shows that this is indeed a spiral attractor and explains the similarity of the shape of the intersection of the numerically obtained attractor with the cross-section and the schematic Figure~\ref{figfoliat4d}b.

\begin{figure}[h!]
\center{\includegraphics[width=0.4\linewidth]{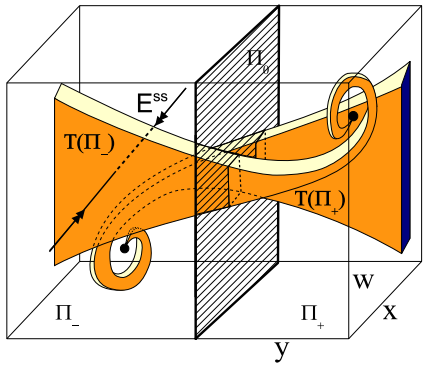} }
\caption{{\footnotesize Schematic model for the Poincar\'e map for a hypothetical case of a pseudohyperbolic spiral attractor with a Lorenz-like geometry.}}
\label{figfol4d_3case2}
\end{figure}

{\bf Remark.}
In the attractor we found in system \eqref{eq:LorenzModified}, the point $M_+$ lies in $\Sigma_+$ and $M_-$ lies in $\Sigma_-$. This is different from what we have for the classical Lorenz attractor (cf. Figs.~\ref{figfoliat4d}b and \ref{figfoliat4d}c). It would be interesting to find examples of pseudohyperbolic spiral attractors for which $M_+\in \Sigma_-$ and $M_+\in\Sigma_+$, like in Fig.~\ref{figfol4d_3case2}. In the kneading diagrams shown in Figs.~\ref{fig:Lorenz4DKneadings} the case $M_+ \in \Sigma_+$ corresponds to region colored in blue, while the Lorenz-like case $M_+ \in \Sigma_-$ corresponds to orange colors (these regions are separated by the curve $l_1$ of a homoclinic butterfly similar to that in the Lorenz model, see Fig.~\ref{fig:foliat}b).

\subsection{Verification of the wild nature of the attractor} \label{sec:kneading}

Our final goal is to demonstrate that the pseudohyperbolic attractor we have found in system \eqref{eq:LorenzModified} for values of parameters close to \eqref{eq:Lor4DParams} is wild, i.e., it admits homoclinic tangencies. The direct search of such tangencies inside the attractor could be a hard computational problem (it requires finding saddle periodic orbits, to construct their invariant manifolds, etc.). Instead, we employ an indirect approach based on the {\it method of kneading diagrams}.

Kneading diagrams were introduced in papers \cite{AShil2012, AShil2014} as a very fast and effective tool for visualization of the complicated bifurcation set corresponding to homoclinic loops to a hyperbolic equilibrium with one-dimensional unstable manifold. We use the kneading diagram to demonstrate the density of parameter values corresponding to the existence of homoclinic loops to the saddle-focus $O$. The latter, by \cite{OvsSh87, OvsSh91}, implies the existence of sought orbits of homoclinic tangencies which pass arbitrarily close to $O$. Since $O$ belongs to the attractor (see Section \ref{sec:geom}), this indicates the wildness of the attractor.

We construct the kneading diagram in the following way. Given a parameter value, we take one of the unstable separatricies of $O$ and use it to build the kneading sequence $s_0,s_1,s_2, \dots$ (by the symmetry, computations with the other separatrix will give equivalent results). If, on this separatrix, the first point corresponding to the maximum of $|x|$ has $x>0$, then we assign $s_1=1$, and if the first maximum of $|x|$ corresponds to $x<0$, then we put $s_0=0$. Repeating the procedure, we can compute the numbers $s_j$ equal to $0$ or $1$ for $j=1,\dots, q$, where $q$ is any aforehand given integer. We always take the right separatrix, so $s_0=1$ and we take it out of the kneading sequence.

\begin{figure}
\center{\includegraphics[width=0.8\linewidth]{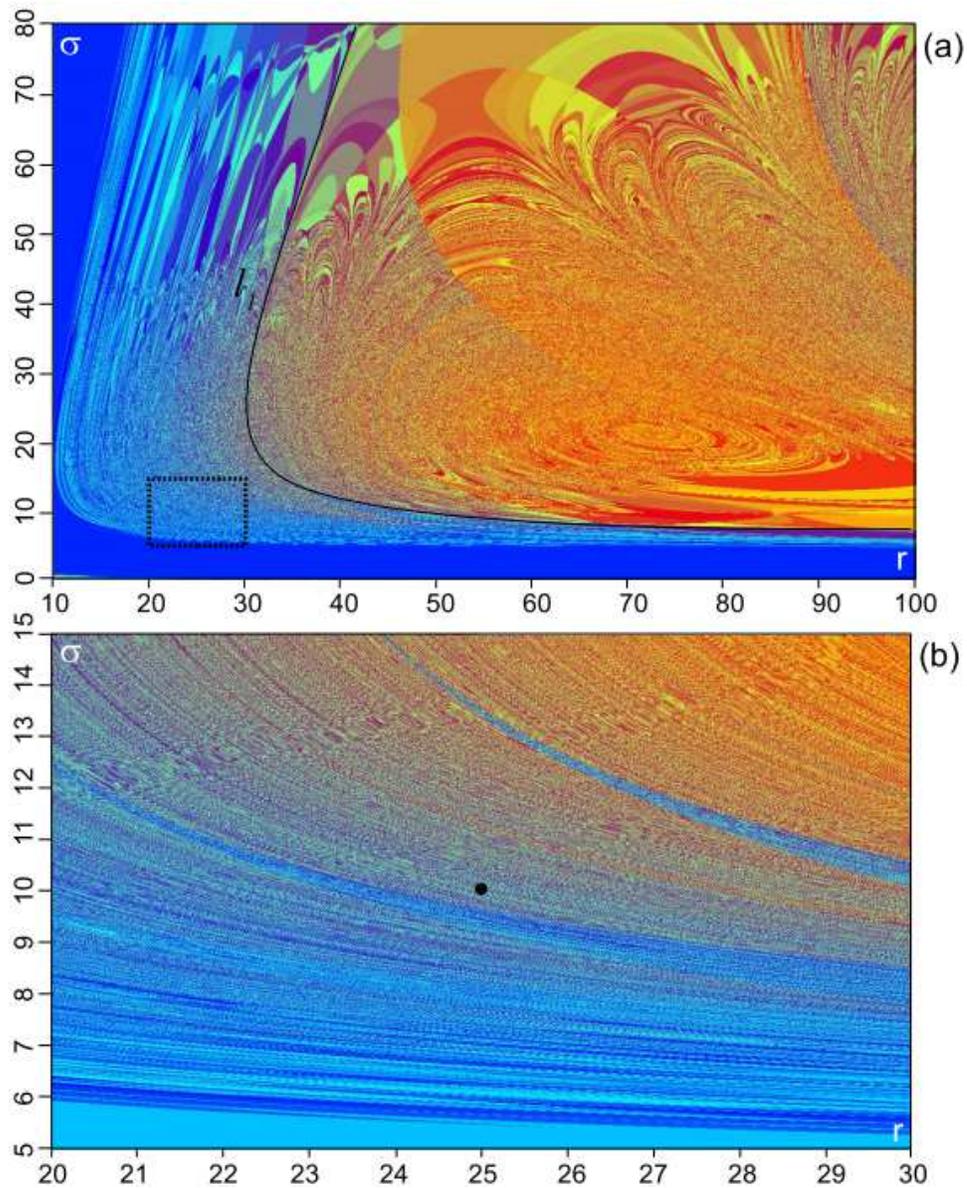} }
\caption{{\footnotesize (a) Kneading diagram for the four-dimensional Lorenz system \eqref{eq:LorenzModified} in the $(r, \sigma)$-plane for fixed $b = 8/3$ and $\mu = 7$; (b) zoomed fragment near the point $(r = 25, \sigma = 10)$. The figure suggests that homoclinic loops to the saddle-focus exist for a dense set of parameter values. The irregular structure of the kneading diagram supports the claim of wildness of the attractor.}}
\label{fig:Lorenz4DKneadings}
\end{figure}

We have shown in Section \ref{sec:geom} that system \eqref{eq:LorenzModified} for values of parameters close to \eqref{eq:Lor4DParams} has a cross-section $\Sigma$ that consists of two disjoint boxes $\Sigma^+$ and $\Sigma^-$. In this region of parameter values, $s_j=1$ means that the $(j+1)$-th point of intersection of the separatrix with $\Sigma$ lies in $\Sigma^+$, and $s_j=0$ means that this point lies in $\Sigma^-$. If for two close parameter values the value of $s_j$ changes while $s_k$ with $k<j$ stay the same, this means that there is a parameter value inbetween which corresponds to the existence of a $j$-round homoclinic loop -- it makes exactly $j$ intersections with $\Sigma$ before closing up\footnote{Note that the existence of a cross-section is important for making such conclusion -- without this the change in the kneading sequence can happen due to events other than formation of a homoclinic loop.}.

For each kneading segment $(s_1, s_2, \dots, s_q)$ we define
$D = \sum_{i=1}^q s_i 2^{q-i}$. Note that $D$ can run integer numbers from $0$ to $2^{q} - 1$, and two length-$q$ kneading segments are equal if and only if the corresponding $D$ values are the same. As we just explained, this means that the boundaries in the parameter space between regions with different values of $D$ correspond to homoclinic loops.
To visualize these boundaries we paint the regions of
different $D$ in different colors -- the resulting picture is the kneading diagram. To do this, we convert each integer from $[0, (2^q - 1)/2)$ to RGB colors following the scheme proposed in \cite{AShil2014}. The values of $D$ from the segment $[0, (2^q - 1)/2)$ are converted to the intensities of red channel, while the blue channel has intensity $0$. The values of $D \in [(2^q - 1)/2, 2^q - 1)]$ are converted to the intensities of the blue channel, while the red channel has intensity $0$. In both these cases the intensity of the green channel takes a random value. This scheme allows to obtain a nicely contrasted picture; we are grateful to Andrey Shilnikov for explaining us these important technical details.

In Fig.~\ref{fig:Lorenz4DKneadings}, kneading diagrams are presented for system \eqref{eq:LorenzModified} in the $(r,\sigma)$ parameter plane for $b=8/3, \mu=7$. Fig.~\ref{fig:Lorenz4DKneadings}a gives a panoramic picture and Fig.~\ref{fig:Lorenz4DKneadings}b shows a zoomed fragment around the point $(r=25,\sigma=10)$. The rapid change of colors
in this figure supports our claim that parameter values corresponding to homoclinic loops to the saddle-focus $O$ are dense. As we explained before, this indicates the presence of homoclinic tangencies inside the attractor.

We also mentioned that bifurcations of homoclinic tangencies cannot be
described by a finite-parameter analysis, meaning that for any finite-parameter unfolding the structure of the bifurcation set is sensitive to small perturbations of the unfolding \cite{GST91, GST93a}. In other words, for any finite-parameter unfolding, the bifurcation set for a system with a wild attractor must have an irregular structure, which is quite convincingly confirmed by the ``blurred'' kneading diagram of Fig.~\ref{fig:Lorenz4DKneadings}b.

\begin{figure}
\center{\includegraphics[width=1.0\linewidth]{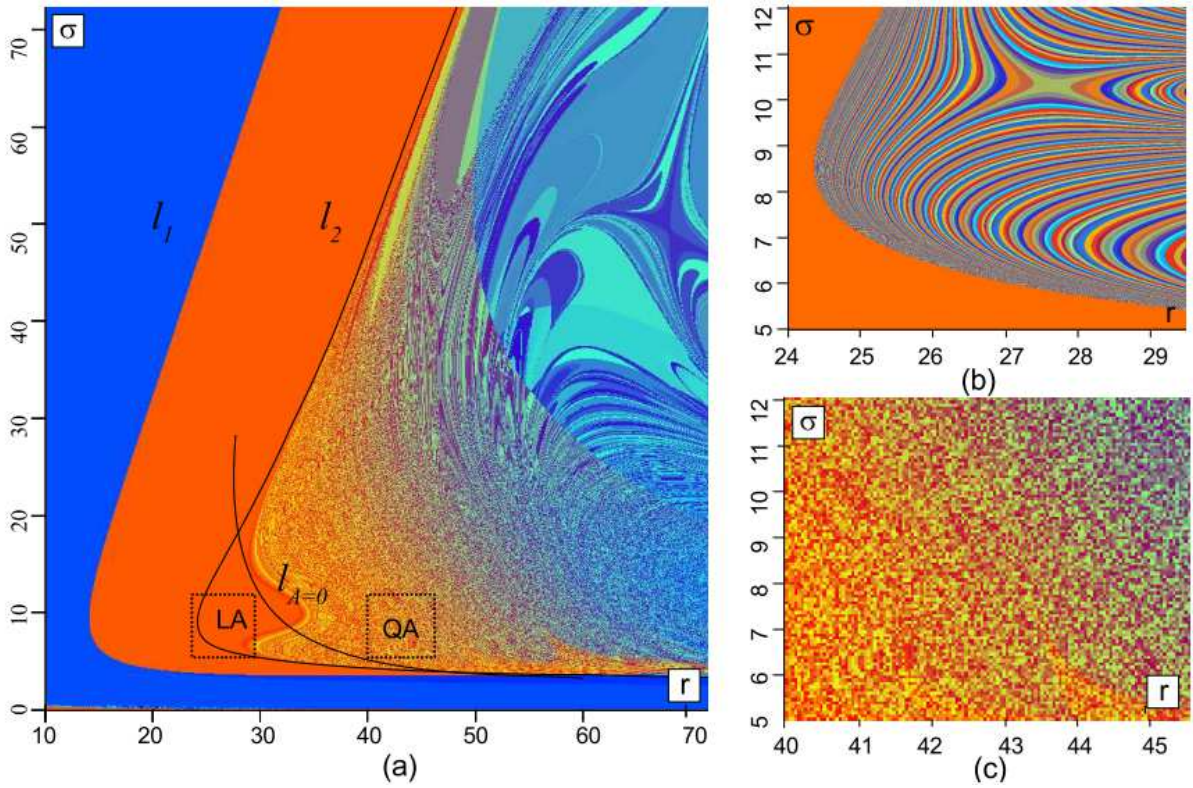} }
\caption{{\footnotesize (a) Kneading diagrams for the Lorenz system \eqref{eq:Lorenz} in the $(r,\sigma)$-plane for $b=8/3$, cf. Fig.~\ref{fig:foliat}b; (b) its zoomed fragment LA near the point $(r=28, \sigma=10)$ where the classical Lorenz attractor exists; (c) zoomed fragment QA near the point $(r=43, \sigma=10)$ -- in this region the system has a quasiattractor.}}
\label{fig:LorenzKneadings}
\end{figure}

In order to illustrate this, we show in Fig.~\ref{fig:LorenzKneadings} the kneading diagram for the classical Lorenz system \eqref{eq:Lorenz} on the $(r, \sigma)$-parameter plane at $b=8/3$ \cite{AShil2012}. In Fig.~\ref{fig:LorenzKneadings}a the diagram of kneading segments of length $q=16$ is presented. We can see that this diagrams is not informative in the domain LA where the Lorenz attractor exists. A more detailed diagram (corresponding to longer kneading sequences) is shown for this domain in Fig.~\ref{fig:LorenzKneadings}b. As we see, the strips with the same kneading segments have a regular structure and are separated from each other by smooth curves corresponding to homoclinic loops to the saddle $O$. It is due to the fact that the kneading sequence is the topological invariant for the Lorenz attractor \cite{Mal85, Mal03}. The situation is drastically changed beyond the curve $l_{A=0}$ where the Lorenz attractor becomes a quasiattractor (see Sec.~\ref{sec:Lor3D}). Kneading diagram becomes here blurred quite similar to what we observe in Fig.~\ref{fig:Lorenz4DKneadings}b, reflecting the fact that homoclinic tangencies appear, see Fig.~\ref{fig:LorenzKneadings}c.

\subsection*{Acknowledgement}

The authors are grateful to A.L. Shilnikov for useful discussions. The authors acknowledge the support of grant No. 075-15-2019-1931 by Russian Ministry of Science and Higher Education, RBBR grants 18-31-20052 and 19-01-00607 (development of geometrical models). DT acknowledges the support by grant EP/P026001/1 by the EPSRC. The theoretical aspects of this work (the theory of pseudohyperbolicity) were carried out with the support by grant RSF No. 19-11-00280. The computational aspects (the development of numerical methods of the verification of pseudohyperbolicity, computation of kneading diagrams) were done with the support by grant RSF No. 19-71-10048.

\end{document}